\documentclass{sapm}

\usepackage{amsmath,amssymb}
\usepackage{graphicx,color}
\usepackage[dvipsnames]{xcolor}
\usepackage{array}
\usepackage{psfrag}
\usepackage{ifthen}
\usepackage{dashrule}

\usepackage{vmargin}
\setpapersize{USletter}
\setmarginsrb{2.5cm}{2.4cm}{2.5cm}{2.4cm}{0.5cm}{0.5cm}{0cm}{0.5cm}

\usepackage{tikz}
\usetikzlibrary{arrows,shapes,backgrounds}

\newcommand{\en}[1]{(\ref{eq:#1})}
\newcommand{\qvect}[1]{\boldsymbol{#1}}
\newcommand{\defeq}{:=}
\newcommand{\order}[1]{{\cal O}\left(#1\right)}
\newcommand{\etal}{{et al.}}

\newcommand{\pfrac}[2]{{\displaystyle \frac{\partial #1}{\partial #2}}}
\newcommand{\half}{\frac{1}{2}}

\newcommand{\Khydraul}{K}

\newcommand{\Dzero}{D_o}
\newcommand{\ystar}{y_{\ast}}

\newcommand{\xstar}{x_{\ast}}
\newcommand{\rstar}{r_{\ast}}
\newcommand{\rstars}{r_{\ast\ast}}

\newcommand{\dtheta}{\Delta\theta}
\newcommand{\eps}{\Theta_{\infty}}
\newcommand{\uinf}{u_{\infty}}
\newcommand{\abstol}{{\tt ABSTOL}}
\newcommand{\reltol}{{\tt RELTOL}}

\newcommand{\tol}{{\tt TOL}}
\newcommand{\Ei}{\operatorname{Ei}}
\newcommand{\li}{\operatorname{li}}

\newcommand{\Matlab}{{\sc Matlab}}
\newcommand{\eulergam}{\Gamma_e}

\newcommand{\vinn}{u}
\newcommand{\vinnz}{\vinn_0}
\newcommand{\vinno}{\vinn_1}
\newcommand{\vinnt}{\vinn_2}
\newcommand{\vinf}{v_{\infty}}

\newcommand{\sstar}{s_{\ast}}

\newcommand{\bbar}{\bar{\beta}}
\newcommand{\sbar}{\bar{\theta}}
\newcommand{\deriv}[2]{{#1}{}^{\prime}}
\newcommand{\dderiv}[2]{{#1}{}^{\prime\prime}}
\newcommand{\ddderiv}[2]{{#1}{}^{\prime\prime\prime}}

\newcommand{\enprime}[1]{(\ref{eq:#1}$^\prime$)}

\definecolor{MyBlue}{rgb}{0,0,1}
\definecolor{MyRed}{rgb}{1,0,0}

\newcommand{\myrevision}[1]{{\color{red}#1}}
\renewcommand{\myrevision}[1]{#1}

\jname{Stud. Appl. Math.}
\jvol{AA}
\jyear{2015}
\doi{10.1146/((please add article doi))}

\begin{document}

\title{Multi-layer asymptotic solution for wetting fronts in porous
  media\ \break with exponential moisture diffusivity}

\author{
  Christopher J.\ Budd\footnotemark[2]
  \and 
  John M.\ Stockie\footnotemark[3]
}

\author[C.\ J.\ Budd and J.\ M.\ Stockie]{Christopher J. Budd and John
  M. Stockie\thanks{Address for correspondence: John M. Stockie,
    Department of Mathematics, Simon Fraser University, 8888 University
    Drive, Burnaby, BC, V5A 1S6, Canada; email: stockie@math.sfu.ca}}
\affil{Department of Mathematical Sciences, University of Bath, Bath, United Kingdom, BA2 7AY} 
\affil{Department of Mathematics, Simon Fraser University, Burnaby, BC, Canada, V5A 1S6}

\maketitle

\begin{abstract}
  We study the asymptotic behaviour of sharp front solutions arising
  from the nonlinear diffusion equation $\theta_t =
  (D(\theta)\theta_x)_x$, where the diffusivity is an exponential
  function $D(\theta)=\Dzero\exp(\beta\theta)$.  This problem arises
  \myrevision{for example} in the study of unsaturated flow in porous
  media where $\theta$ represents the liquid saturation.  For 
  physical parameters corresponding to actual porous media, the
  diffusivity at the residual saturation is $D(0)=\Dzero\ll 1$ so that
  the diffusion problem is nearly degenerate.  Such problems are
  characterised by wetting fronts that sharply delineate regions of
  saturated and unsaturated flow, and that propagate with a well-defined
  speed.  Using matched asymptotic expansions in the limit of large
  $\beta$, we derive an analytical description of the solution that is
  uniformly valid throughout the wetting front.  This is in contrast
  with most other related analyses that instead truncate the solution at
  some specific wetting front location, which is then calculated as part
  of the solution, and beyond that location the solution is undefined.
  Our asymptotic analysis demonstrates that the solution has a
  four-layer structure, and by matching through the adjacent layers we
  obtain an estimate of the wetting front location in terms of the
  material parameters describing the porous medium.  Using numerical
  simulations of the original nonlinear diffusion equation, we
  demonstrate that the first few terms in our series solution provide
  approximations of physical quantities such as wetting front location
  and speed of propagation that are more accurate (over a wide range of
  admissible $\beta$ values) than other asymptotic approximations
  reported in the literature.
\end{abstract}




\section{Introduction}
\label{sec:intro}

\begin{figure}[tbhp]
  \begin{center}
    \leavevmode
    \includegraphics[width=0.5\textwidth,clip]{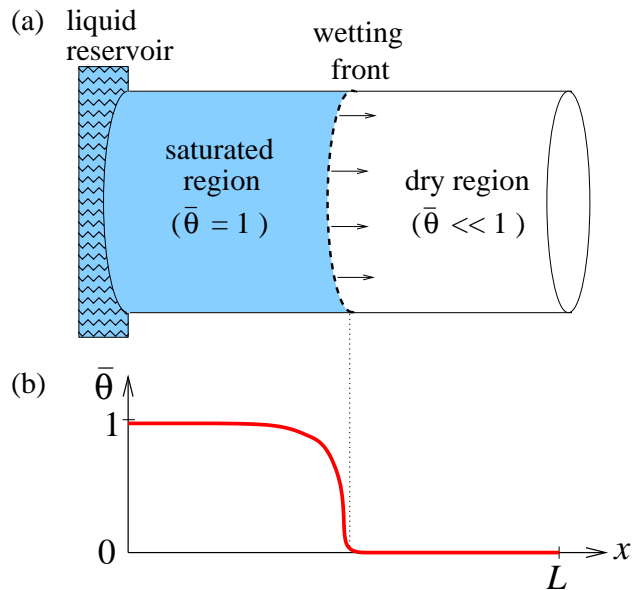}
    \caption{\myrevision{(a) A cylindrical porous medium with the left
        end immersed in a water reservoir, depicting the progress of a
        wetting front to the right owing to capillary action.  (b) Plot
        of rescaled saturation $\sbar = (\theta-\theta_o) /
        (\theta_i-\theta_o)$ along the length of the cylinder from $x=0$
        to $L$.}}
    \label{fig:geometry}
  \end{center}
\end{figure}

Problems with exponential diffusivity arise commonly in the study of
water transport in variably-saturated porous media such as soil, rock or
building materials~\cite{brutsaert-1979}.  An example of such,
\myrevision{expressed in dimensionless form,} is the one-dimensional
nonlinear diffusion problem
\begin{subequations}\label{eq:diffusion}
  \begin{gather}
    \pfrac{\theta}{t} = \pfrac{}{x} \left( D(\theta) \pfrac{\theta}{x}
    \right),
    \label{eq:diffusion-eqn} \\
    \theta(0,t)=\theta_i \quad \text{and} \quad
    \theta(L,t)=\theta_o \quad \text{for $t\geqslant 0$}, \label{eq:diffusion-bcs} \\
    \theta(x,0)=\theta_o \quad \text{for $0<x<L$}, \label{eq:diffusion-ic} 
  \end{gather}%
  %
  %
  where the solution $\theta(x,t)$ is called \emph{saturation} and
  represents the volume fraction of pore space occupied by liquid.  We
  are concerned here with the case when \myrevision{the dimensionless
    variable} $L\gg 1$ and the diffusivity is an {\em exponential
    function} of the solution having the form
  \begin{gather}
    D(\theta)=\Dzero e^{\beta\theta},
    \label{eq:diffusion-Dexp}
  \end{gather}%
  %
  %
\end{subequations}
where $\Dzero$ and $\beta$ are constants satisfying $0<\Dzero\ll 1$ and
$\beta\gg 1$.  \myrevision{This is a reasonable model for horizontal
  infiltration problems such as that pictured in
  Figure~\ref{fig:geometry}a, where liquid is taken up by capillary
  action within a dry porous medium and gravity can be neglected.}  Many
experimental studies of porous media have been performed in which an
exponential diffusion ansatz provides a good fit with measured data,
mostly in the context of water transport in soil and
rock~\cite{brutsaert-1979, clothier-white-1981, miller-bresler-1977,
  reichardt-nielsen-biggar-1972}, but also for other porous materials
such as wood~\cite{simpson-1993}, brick~\cite{pel-1995} or
concrete~\cite{leech-lockington-dux-2003}.  Exponential diffusion has
also been identified in other transport phenomena as diverse as
\myrevision{heat conduction~\cite{cooper-1971}, 
  optical lithography~\cite{wagner-1950}, 
  solvent transport in polymers~\cite{hansen-1980}, and diffusion of
  impurities in oxides~\cite{wei-wuensch-1976}}.  The significance of
this type of diffusion coefficient is that if $\theta_i > \theta_o$ and
if $\beta$ is even moderately large then $D(\theta)$ varies rapidly with
$\theta$, thereby causing solutions of \en{diffusion-eqn} to develop a
steep \myrevision{interface, called a {\em wetting front} in the context
  of porous media flow,} that is associated with localised high
curvature and an exponential change in the solution gradient
\myrevision{(see Figure~\ref{fig:geometry}b)}.  The aim of this paper is
to perform an asymptotic analysis of this phenomenon and in particular
to study self-similar solutions of \en{diffusion-eqn} using a
multi-layer asymptotic expansion, supported by numerical
calculations. The asymptotic theory is especially subtle owing to the
exponential change in the solution gradient and yields sharp estimates
that agree well with numerical simulations.

\subsection{Overview of Previous Work}

The nonlinear diffusion equation~\en{diffusion-eqn} has been studied in
considerable detail for diffusivities $D(\theta)$ having a variety of
functional forms.  Particular emphasis has been placed on the special
case of a power law, $D(\theta)=\myrevision{\Dzero}\theta^m$ with $m>0$,
for which \en{diffusion-eqn} is called the \emph{porous medium equation}
or PME; an extensive literature exists for the PME that is thoroughly
covered in the review by V\'{a}zquez~\cite{vazquez-2007}. In a typical
wetting scenario the PME is supplemented by boundary conditions
$\theta(0,t)=1$ and $\theta(\infty,t)=0$, in which case the solution for
a power-law diffusivity is well known to have compact support and to
consist of a front propagating to the right with speed proportional to
$t^{-1/2}$.  Ahead of the front, the solution is identically zero and
when $m\geqslant 1$ there is a discontinuity in the first derivative
$\theta_x$ at the point where the front meets the $x$-axis; otherwise,
when $0<m<1$, the wetting front meets the leading edge solution
smoothly.  The trailing edge of the front, on the other hand, always
experiences a smooth transition similar to that shown in the saturation
profile in Figure~\ref{fig:geometry}b. Analytical results have been
derived for other types of diffusion coefficient, for example by
Shampine~\cite{shampine-1973b} who provides an existence-uniqueness
proof for a general class of diffusivity functions with $D(\theta)$
required to be a continuously differentiable function on $0\leqslant
\theta\leqslant 1$, satisfying $D(0)=0$ and $D(\theta)>0$ when
$\theta>0$.
This and other analytical studies are characterised by the fact that
they pertain to \emph{degenerate diffusion} for which the diffusivity
vanishes identically at zero saturation.

\begin{figure}[tbhp]
  \centering
  \psfrag{theta}[cc][bc]{$\theta$}
  \psfrag{D\(theta\)}[bc][cc]{$D(\theta)/\Dzero$}
  \psfrag{y}[cc][bc]{$x$}
  \psfrag{theta\(y\)}[bc][cc]{$\theta$}
  \begin{tabular}{ccc}
    (a) Rescaled diffusion coefficients & \qquad & 
    (b) Saturation profiles \\
    \includegraphics[width=0.35\textwidth,clip]{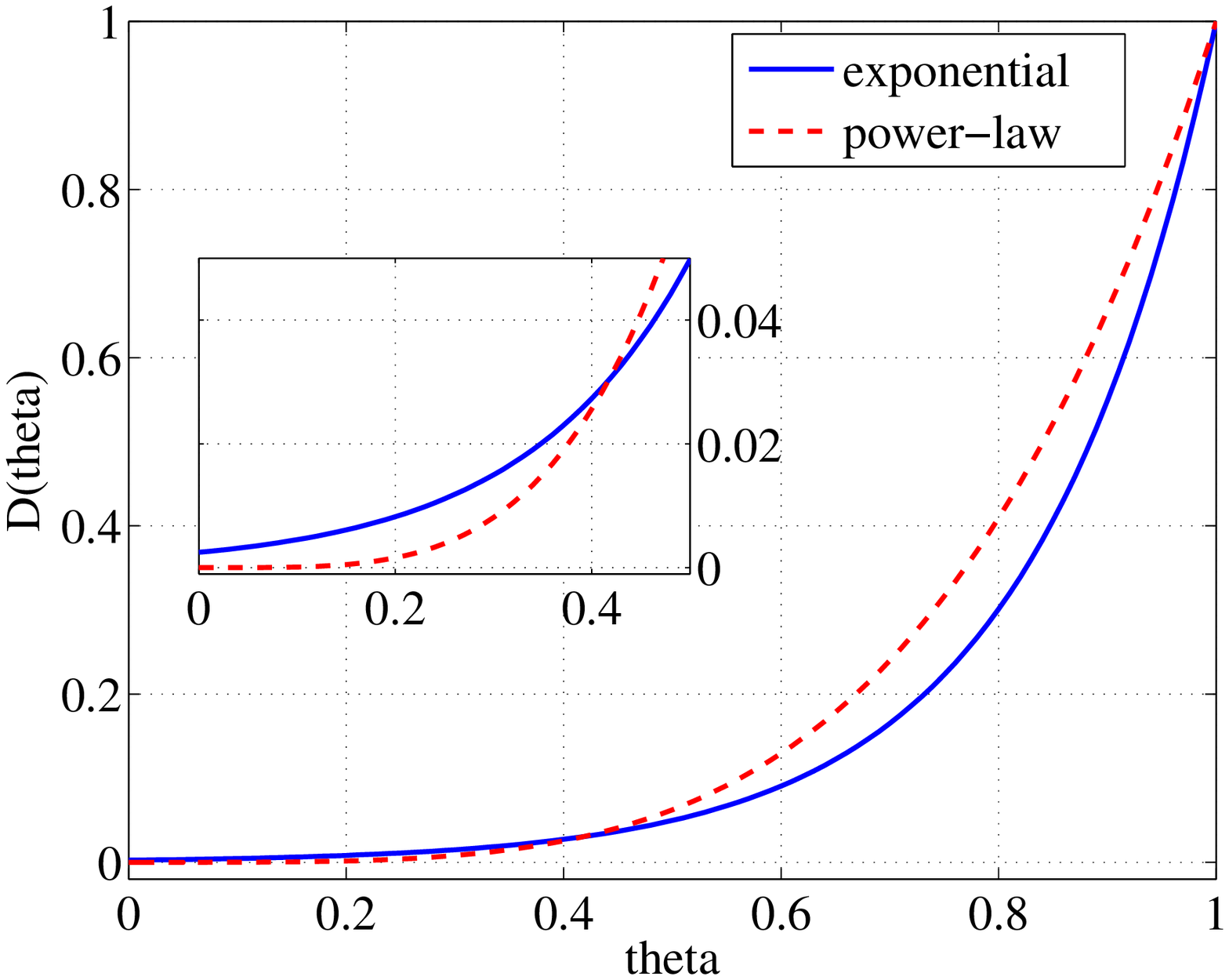}
    & & 
    \begin{tikzpicture}
      \node (corner-exp) at (,) {\Large Note};
      \node (corner-power) at (,) {\Large Water};
      %
      \node[anchor=south west,inner sep=0] (image) at (0,0)
      {\includegraphics[width=0.35\textwidth,clip]{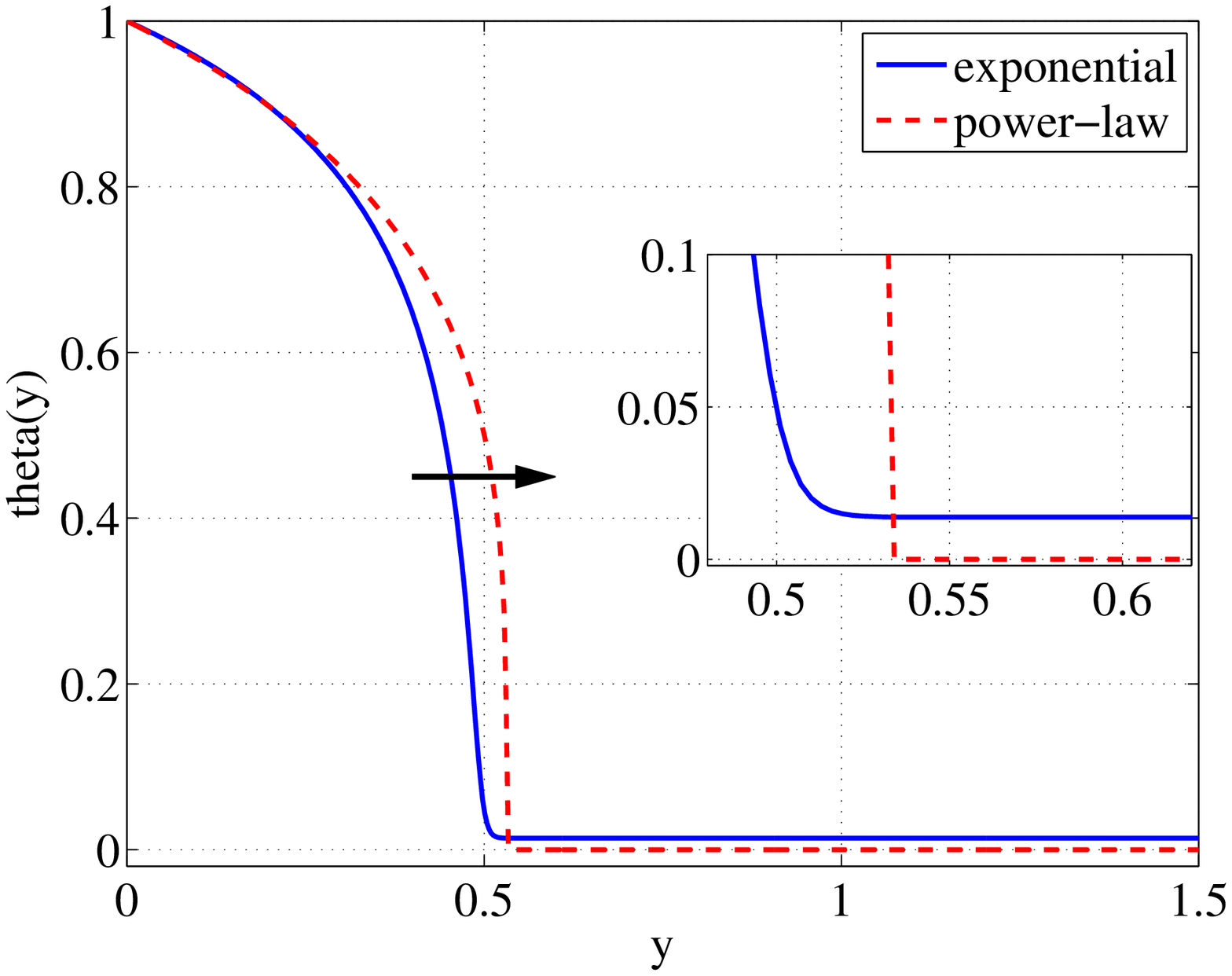}};
      \begin{scope}[x={(image.south east)},y={(image.north west)}]
        \node[ForestGreen] (wet-front1) at (0.87,0.65) {\scriptsize Wetting};
        \node[ForestGreen] (wet-front2) at (0.90,0.6) {\scriptsize front};
        \node[ForestGreen] (wet-front3) at (0.91,0.54) {\scriptsize $\xstar$};
        \draw[-latex,ForestGreen,line width=1pt] (wet-front2) to[out=180, in=40] (0.675,0.485);
        \draw[-latex,ForestGreen,line width=1pt] (wet-front2) to[out=180, in=35] (0.735,0.43);
        \draw[ForestGreen] (0.665,0.475) circle (1.5pt);
        \draw[ForestGreen] (0.725,0.420) circle (1.5pt);
      \end{scope}
    \end{tikzpicture}
  \end{tabular}
  \caption{Comparison of solutions to the nonlinear diffusion problem
    for both power-law and exponential diffusion coefficients,
    $D(\theta)=\Dzero \theta^4$ and $D(\theta)=\Dzero
    \myrevision{e^{6\theta}}$.  (a) On the left is the rescaled
    diffusivity $D(\theta)/\Dzero$.  (b) On the right are the
    corresponding saturation profiles $\theta(x,t)$ computed numerically
    at some fixed time $t>0$.  \myrevision{The black arrow indicates the
      direction of the wetting front motion, from left to right.  The
      green arrows within the inset plot denote the wetting front
      $\xstar$ for each solution, located at the point of maximum
      curvature.}}
  \label{fig:power}
\end{figure}

In contrast with this previous work, the diffusion coefficient we
consider in this paper is an exponential function of saturation that,
although small, is still bounded away from zero; hence the solution
remains everywhere classical and all disturbances propagate with
infinite speed, analogous to solutions of the ``usual'' heat equation.
When the PME wetting scenario described above is repeated for an
exponential diffusivity, the solution remains non-zero for all times
$t>0$, even when the initial conditions have bounded support.
Nonetheless, the problem can be \emph{nearly degenerate} in the sense
that $D(0)\ll 1$ whereas $D(\theta)=\order{1}$ for $\theta$ away from
zero.  For example, typical parameter values for a porous soil have been
estimated as $\Dzero\approx 2\times 10^{-9}\;m^2/s$ and $\beta\approx
20$~\cite{reichardt-nielsen-biggar-1972}.\ \
For such diffusion coefficients, the solution inherits features that are
qualitatively similar to the degenerate problem, most notably a steep
wetting front that propagates with finite speed and that is associated
with localised high curvature (see Figure~\ref{fig:power}(b)).  We
define the wetting front location $\xstar(t)$ to be the point where the
curvature of $\theta(x,t)$ is greatest, and for $x > \xstar$ we have
that $\theta_x \ll 1$.  One feature that distinguishes the exponential
diffusion equation from the PME is the fact that $\theta$ exhibits more
rapid variation than for a power law, which in turn has a substantial
impact on the shape of the solution close to the wetting front. A direct
comparison is afforded by Leech et al.'s study of
concrete~\cite{leech-lockington-dux-2003} wherein they use experimental
data to fit both power-law and exponential diffusivities, yielding
respectively $D(\theta)=\Dzero\theta^4$ and $D(\theta)=\Dzero
\myrevision{e^{6\theta}}$.  In Figure~\ref{fig:power}, we plot the
diffusion coefficients and corresponding saturation profiles for these
two choices of $D(\theta)$.  While the qualitative features of both
solutions are similar, there is a significant difference in both the
steepness and location of the wetting front, as well as the sharpness of
the corner (refer to the zoomed-in region in Figure~\ref{fig:power}(b)).
Another distinguishing feature is that the solution in the exponential
case is characterized by a small nonzero saturation ahead of the wetting
front \myrevision{due to the fact that $D(0)\neq 0$}.

Following up on these experimental studies, several authors have derived
theoretical results for the exponential diffusion problem.  Crank's
book~\cite{crank-1975} provides a comprehensive treatment of analytical
solutions for nonlinear diffusion equations circa 1975 with many forms
of the diffusivity function.  In particular, Crank derives a similarity
solution for the exponential diffusion case (following the work of
Cooper~\cite{cooper-1971}) that reduces the problem to a second-order
ordinary differential equation (ODE) which he then solves numerically --
indeed, this is the same ODE that we will present later in
Section~\ref{sec:self-similar}.  An alternate approach using functional
iteration has been developed based on an integral formulation of the
governing equations by Parslow
\etal~\cite{parslow-lockington-parlange-1988}.

In contrast with these iterative or numerical solution methods, we are
interested here in developing an asymptotic series representation of the
solution.  One study of particular relevance is due to
Babu~\cite{babu-1976a} who derived an asymptotic solution by making use
of the simplifying assumption that both saturation and diffusivity drop
to zero ahead of a certain wetting front location.  Although Babu's
asymptotic estimates are reasonably accurate, we will demonstrate that
his approach of truncating the solution at the wetting front introduces
significant errors that can be reduced by considering an alternate
approach that incorporates the effect of the extremely small but still
nonzero diffusivity values ahead of the front.  Other alternate series
expansions based on an integral form of the exponential diffusion
equation were derived by Parlange and co-workers in
\cite{parlange-1973,parlange-etal-1992}.
Finally, we point out a connection to the work of Elliott
\etal~\cite{elliott-etal-1986} who applied methods from singular
perturbation theory to analyse the detailed structure of the transition
region for the power-law diffusivity \myrevision{$D(\theta) =
  \Dzero\theta^m$} in the limiting case of $m\to\infty$, which is known
as the \emph{mesa problem}.  Although this problem is still degenerate,
the diffusivity in the large-$m$ limit experiences a rate of increase in
$\theta$ approaching that of an exponential function; therefore, the
analysis for the mesa problem can be considered as a prototype for
problems such as \en{diffusion}.  \myrevision{Asymptotic solutions were
  developed for a semi-conductor dopant diffusion problem by
  King~\cite{king-1988} for the two limits $m\to 0$ (constant $D$) and
  $m\to\infty$ (mesa), and King and Please~\cite{king-please-1986}
  applied singular perturbation theory to obtain a matched asymptotic
  solution for the case $m=1$ (linear $D$); however, the asymptotic
  structure of the solution for the exponential diffusion problem has
  not yet been explored in similar detail.}

\subsection{Summary of Main Results}

We will demonstrate in this paper that a multi-layer asymptotic
expansion is capable of yielding an accurate estimate of the wetting
front location, provided that the exponentially varying solution is
properly resolved close to the front. In particular, we will show that
there is a self-similar solution \myrevision{$\Theta(y)$ that can be
  written in terms of the dimensionless similarity variable
  \begin{gather*}
    y = \frac{x}{\sqrt{2tD(\theta_i)}}. 
  \end{gather*}
  The corresponding wetting front location $\ystar$ is a constant and
  has the asymptotic expansion}
\begin{gather*}
  \ystar = \frac{1}{\gamma} + \frac{1}{2 \gamma^3} + \frac{11}{12
    \gamma^5} + \order{\frac{1}{\gamma^7}}. 
\end{gather*}
The large parameter
$\gamma =  - \Theta^\prime(0)$
is related to the physical constant $\beta$ through the asymptotic
expression
\begin{gather*}
  \bbar \equiv \beta \left( \theta_i - \theta_o \right) = \gamma^2
  + \frac{1}{2} + \frac{\alpha_3}{\gamma^2} + \order{\frac{1}{\gamma^4}}
  \gg 1, 
\end{gather*}
where $\alpha_3 \approx \frac{1}{12}$ is obtained numerically.  Our
asymptotic solution is distinct from other approximations derived in the
literature in that it provides insight into the detailed structure of
the wetting front, as well as yielding estimates of quantities such as
$\ystar$ that are accurate over a wide range of parameters.
\myrevision{The estimate of wetting front location could be of
  particular interest to engineers since it expresses an easily
  measurable quantity in terms of physical parameters.}

The remainder of this paper is structured as follows.  In
Section~\ref{sec:background} we motivate the problem under study by
using the example of water transport in unsaturated porous media,
\myrevision{although we stress that these results are equally applicable
  to other nonlinear diffusion problems having an exponential
  diffusivity}.  In Section~\ref{sec:self-similar} we introduce a
similarity transformation that permits us to recast problem as an ODE
initial value problem, for which $\gamma\gg 1$ \myrevision{is related
  to} the magnitude of the initial slope of the similarity solution.
Section~\ref{sec:asymptotic} derives our key results, consisting of a
four-layer asymptotic expansion for the similarity solution, expressed
as a series in $\gamma$ on each layer.  Matching the asymptotic
expressions then yields series approximations for the quantities of
physical interest such as the saturation and curvature at the sharp
corner in the wetting front, as well as the location $\xstar$ of the
front itself.  In Section~\ref{sec:numerical} we perform a detailed
comparison of our asymptotic solution to other approximations in the
literature, as well as validating the results with careful numerical
simulations of the original governing partial differential equation.

\section{Physical Background}
\label{sec:background}


Water transport in a saturated porous medium is well-known to obey
Darcy's law~\cite{bear-1988}, $\qvect{U} = -\Khydraul \nabla \Phi$,
which states simply that the liquid velocity $\qvect{U}$ (in $m/s$) is
proportional to the gradient of total hydraulic potential $\Phi$ (in
$m$).  The proportionality constant $\Khydraul$ is known as the
hydraulic conductivity and has units of $m/s$.  However, many porous
media flows are unsaturated, which means that the pore volume is only
partially filled with liquid, and in this case it is necessary to
introduce the liquid volume fraction or saturation, $\theta$.  The
conductivity in a variably saturated porous medium is typically assumed
to depend on the local saturation, leading to an extended Darcy's law,
$\qvect{U} = -\Khydraul(\theta) \, \nabla \Phi$, that includes a
saturation-dependence in the hydraulic conductivity.  When this
expression for velocity is substituted into the continuity equation
\begin{gather*}
  \pfrac{\theta}{t} + \nabla \cdot \qvect{U} = 0,
\end{gather*}
we obtain an evolution equation for $\theta$
\begin{gather}
  \pfrac{\theta}{t} = \nabla \cdot \left( \Khydraul(\theta) \nabla
  \Phi \right),   
  \label{eq:richards}
\end{gather}
which is known as the \emph{Richards equation}.  If both $\Phi$ and
$\Khydraul$ are assumed to be single-valued functions of saturation,
then we can define $D(\theta) \defeq \Khydraul(\theta) \, (d\Phi/d\theta)$
after which Eq.~\en{richards} reduces to the familiar nonlinear
diffusion equation
\begin{gather*}
  \pfrac{\theta}{t} = \nabla \cdot \left( D(\theta) \nabla \theta \right),  
  \label{eq:diffusion-nd}
\end{gather*}
for which Eq.~\en{diffusion-eqn} is the 1D version.  The function
$D(\theta)$ is referred to as the \emph{moisture diffusivity} and has
units of $m^2/s$.  As mentioned in the Introduction, the exponential
form \en{diffusion-Dexp} for diffusivity is motivated by the study of
certain soils and porous building materials, for which an exponential
function provides a good fit to experimental data.


We consider an idealised, cylindrical geometry depicted in
Figure~\ref{fig:geometry}a that is consistent with the samples typically
used in experimental studies of moisture transport.
The porous cylinder is initially dry and has one end placed in a water
reservoir.  \myrevision{We are interested in problems such as the
  \emph{horizontal infiltration} scenario pictured, or where capillary
  forces dominate over gravity, so that water is absorbed into the
  porous medium by capillary action alone.  Water progresses into the
  sample as a nearly planar {wetting front}, and so we can assume} that
the flow is uni-directional and governed by the one-dimensional
diffusion equation~\en{diffusion-eqn}.  The length of the sample is
denoted by $L$, which is presumed large in relation to the diameter so
that $L\gg 1$.

In order to treat a general class of wetting scenarios we take samples
that are never completely dry, which corresponds to the usual situation
wherein a porous medium undergoes re-wetting after an initial
wetting/draining cycle.  Consequently, there exists a non-zero minimum
or residual saturation $\theta=\theta_o$ deriving from water that is
trapped in isolated portions of the porous matrix and which cannot be
displaced by capillary action.  Furthermore, we introduce a maximum
saturation $\theta_i\leqslant 1$ that cannot be exceeded owing to
micropores that are too small to allow water to penetrate, no matter how
large the capillary force.  Consequently, the saturation $\theta$
satisfies $0 < \theta_o \leqslant\theta\leqslant\theta_i\leqslant 1$,
while the boundary and initial conditions are as specified in
\en{diffusion-bcs} and \en{diffusion-ic}.  Water content is usually
reported in the literature in terms of \emph{reduced saturation}
\begin{gather*}
  \sbar = \frac{\theta-\theta_o}{\theta_i-\theta_o}
  = \frac{\theta-\theta_o}{\dtheta} ,
\end{gather*}%
%
%
for which the boundary and initial conditions reduce to 
\begin{gather}
  \sbar(0,t)=1, \qquad 
  \sbar(L,t)=0 \quad \text{and} \quad  
  \sbar(x,0)=0.  
\end{gather}
A picture of a typical wetting front is displayed in
Figure~\ref{fig:geometry}b in terms of the rescaled saturation
variable $\sbar$.


\section{Self-Similar Solution}
\label{sec:self-similar}

We are now interested in finding a self-similar solution of the
nonlinear diffusion equation on the half-space $x\in [0,\infty]$,
\myrevision{when $L\to\infty$, and which satisfies the Dirichlet
  boundary condition $\theta(0,t) = \theta_o$ constant}.  The problem
\en{diffusion} is invariant under the scaling transformation $y =
x t^{-1/2}$,
which suggests seeking a solution of the form
\myrevision{$\theta(x,t)=\phi(y)$, where the specific form of the
  similarity variable
  \begin{gather}
    y = \frac{x}{\sqrt{2tD(\theta_i)}}  
    \label{eq:boltzmann}
  \end{gather} 
  is chosen to eliminate certain constants from the ODE and
  boundary conditions.}
We then define the new dependent variable 
\begin{gather}
  \myrevision{\Theta \defeq e^{\beta (\phi(y) - \theta_i)}}
  \qquad \text{or equivalently} \qquad
  \Theta \defeq e^{\bbar (\sbar-1)} ,
  \label{eq:similar-sat}
\end{gather}
where 
\begin{gather*}
\bbar \defeq \beta\dtheta = \beta(\theta_i-\theta_o).
\end{gather*}
This form of solution is consistent with the second boundary
condition~\en{diffusion-bcs} \myrevision{provided that
  $L/\sqrt{D(\theta_i)t}\gg 1$,} which motivates our taking $L\to
\infty$ in what follows.
%
%
This simplifies the problem further by eliminating the exponential
diffusion coefficient \en{diffusion-Dexp}, leading to the following ODE
boundary value problem for $\Theta(y)$:
\begin{subequations}\label{eq:usim}
  \begin{align}
    \Theta \dderiv{\Theta}{y} &= -y \deriv{\Theta}{y}
    \qquad\text{for $0<y<\infty$}, \label{eq:usim-u}\\
    \Theta(0) &= 1,\label{eq:usim-ic1}\\
    \Theta(\infty) &= \eps \defeq e^{-\bbar} \label{eq:usim-bcL}.
  \end{align}%
\end{subequations}%
%
%
For the porous media of interest here, $\bbar$ is significantly greater
than one so that the parameter $\eps$ satisfies $0<\eps\ll 1$.

Rather than solving Eqs.~\en{usim} directly it is convenient, for both
the numerical and the asymptotic calculations, to reformulate this
boundary value problem as an initial value problem, in which the
boundary condition \en{usim-bcL} is replaced by a second initial
condition of Neumann type given by
\begin{gather}
  \deriv{\Theta}{y}(0) 
  \;\myrevision{= \beta\phi^\prime(0)}
  = -\gamma, \tag{\ref{eq:usim-bcL}${}^{\prime}$}
  \label{eq:usim-ic2}
\end{gather}
where \myrevision{$-\gamma$ represents the initial slope which we will
  assume to be large}.\ \
%
%
From this point on, we use \enprime{usim} to refer to the initial value
problem consisting of equations \en{usim-u}, \en{usim-ic1} and
\en{usim-ic2}.  The asymptotic analysis performed in this paper
considers the limit of $\gamma \to \infty$ which we will see shortly is
equivalent to taking $\eps \to 0$.  The quantity $\gamma$ is not known
\emph{a priori} in terms of the physical parameters, and so one of the
primary results of this paper will be to first derive asymptotic
expressions for $\bbar$ and $\ystar$ in terms of $\gamma$ and to then
invert them to give $\gamma$ in terms of the physical quantities.

\myrevision{Before moving on, we address the suitability of applying
  alternate boundary conditions of Neumann (or flux) type.  First of
  all, it follows from equations \en{usim-u} and \en{usim-ic1}, along
  with a later result indicating $\Theta^{\prime\prime}(\infty)\sim 0$,
  that the Dirichlet problem above has an implicit no-flux (zero
  Neumann) condition at infinity.  Hence, there is no need to consider
  this case separately.  When the left boundary condition is replaced
  with a flux contition of the form $\theta_x(0,t)=a$, no similarity
  solution is permitted unless $a=0$, and in the zero-flux case there is
  only the trivial solution $\Theta\equiv\eps$.  Therefore, it is
  sufficient to focus on the Dirichlet problem in \en{usim}.}

\subsection{Preliminary Numerical Simulations}

\begin{figure}[tbhp]
  \psfrag{y}[cc][bc]{$y$}
  \psfrag{u\(y\)=exp\(B*\(theta-thetaI\)\)}[bl][cc]{$\Theta(y)$}
  \psfrag{x}[cc][bc]{$x$}
  \psfrag{theta\(x,t\)}[bc][cc]{$\sbar(x,t)$}
  \begin{center}
    \leavevmode
    \begin{tabular}{ccc}
      (a)\ \ $\bbar=2$ & (b)\ \ $\bbar=4$ & (c)\ \ $\bbar=8$ \\
      $(\gamma=1.166,\;\eps=0.135)$  &  
      $(\gamma=1.850,\;\eps=0.0183)$ & 
      $(\gamma=2.736,\;\eps=0.000335)$ \\
      \includegraphics[width=0.30\textwidth,clip]{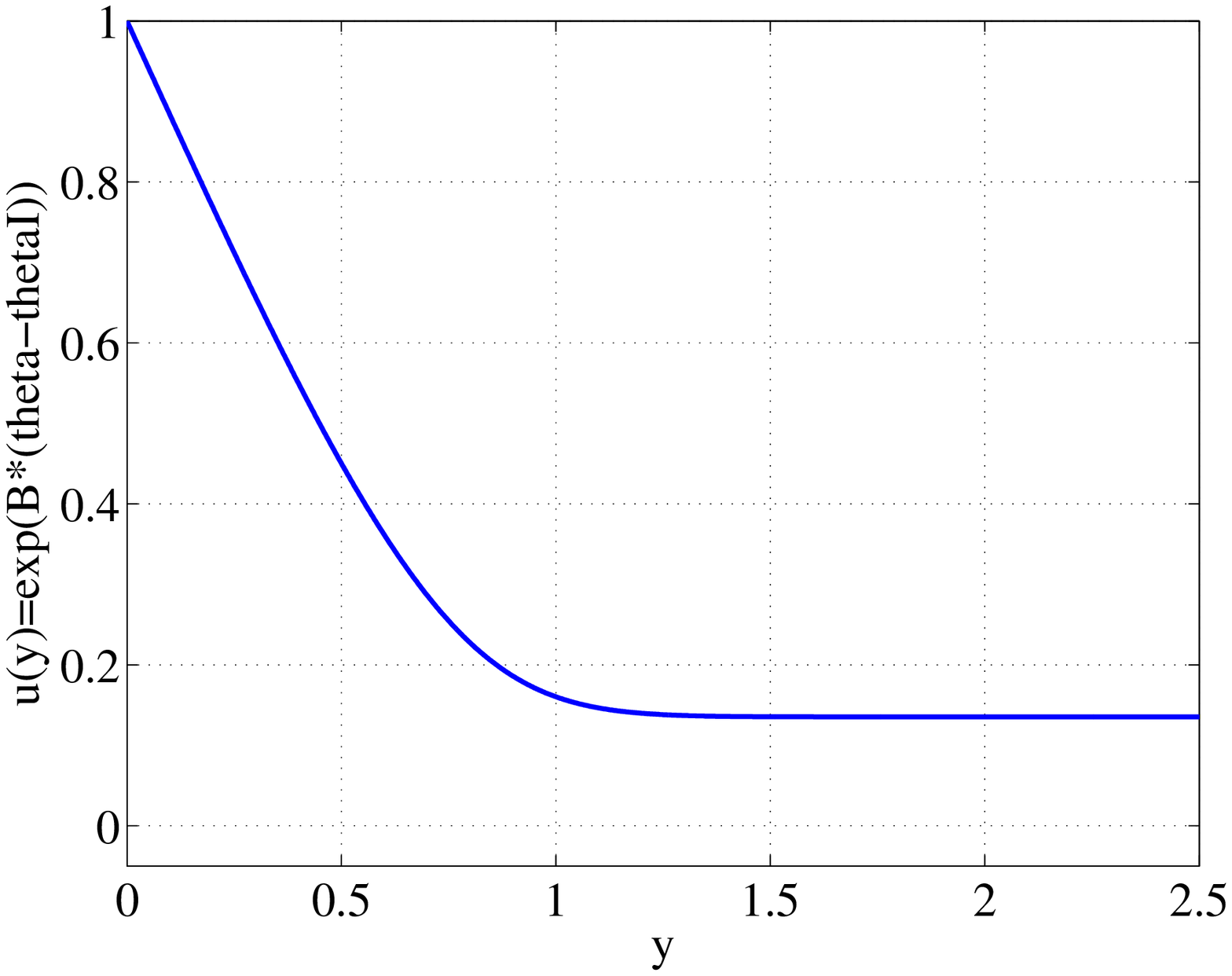} &
      \includegraphics[width=0.30\textwidth,clip]{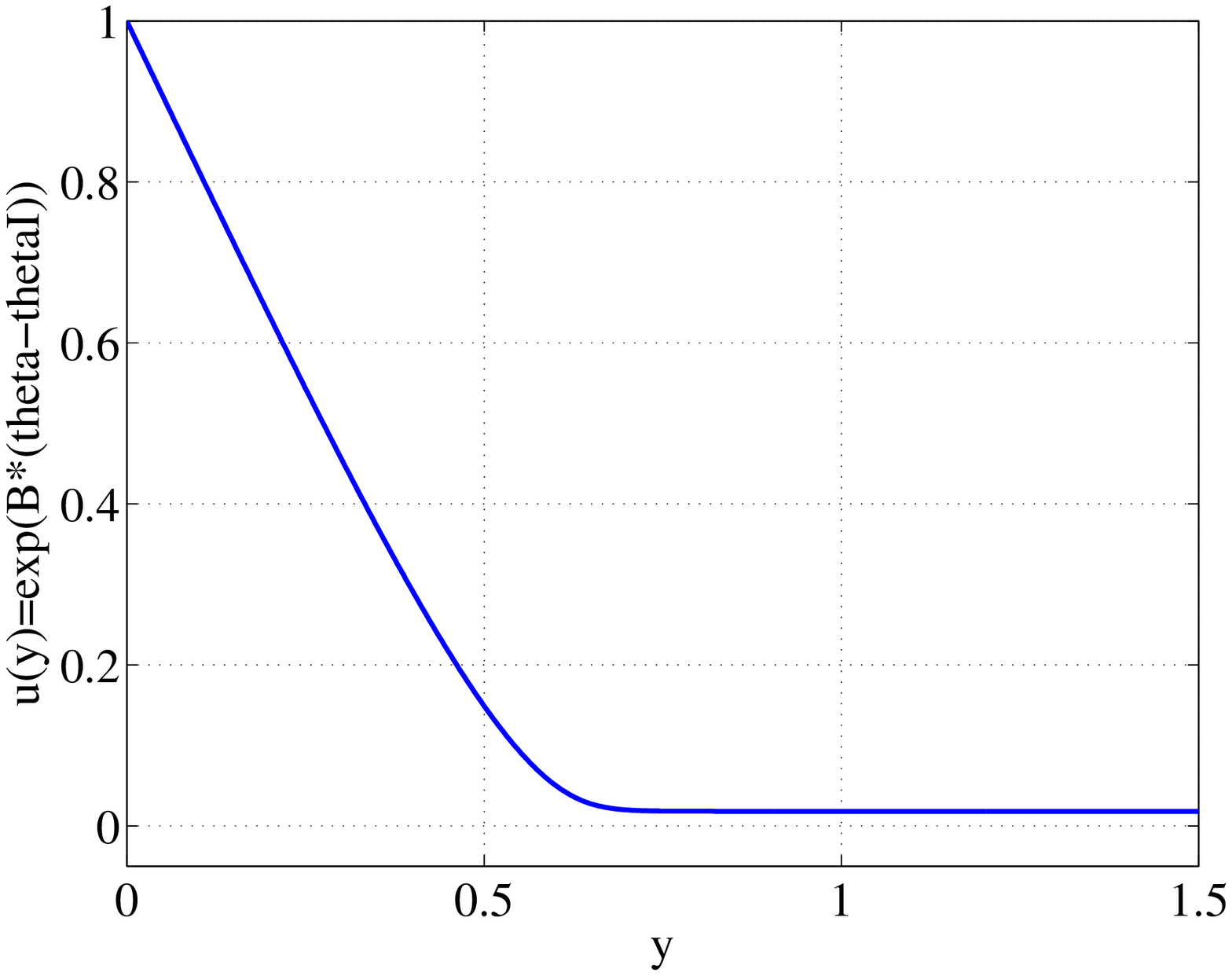} &
      \includegraphics[width=0.30\textwidth,clip]{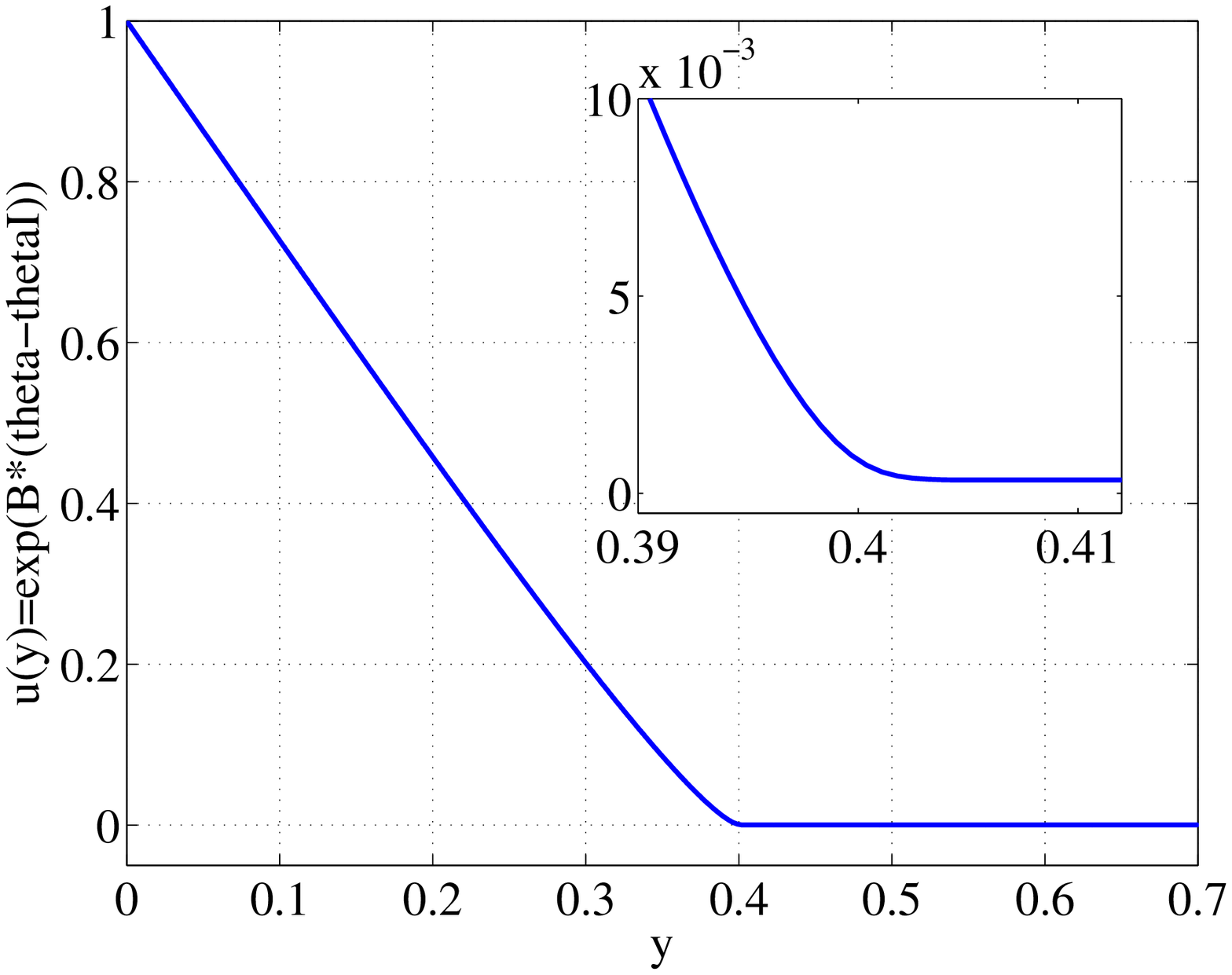}
      \\
      \includegraphics[width=0.30\textwidth,clip]{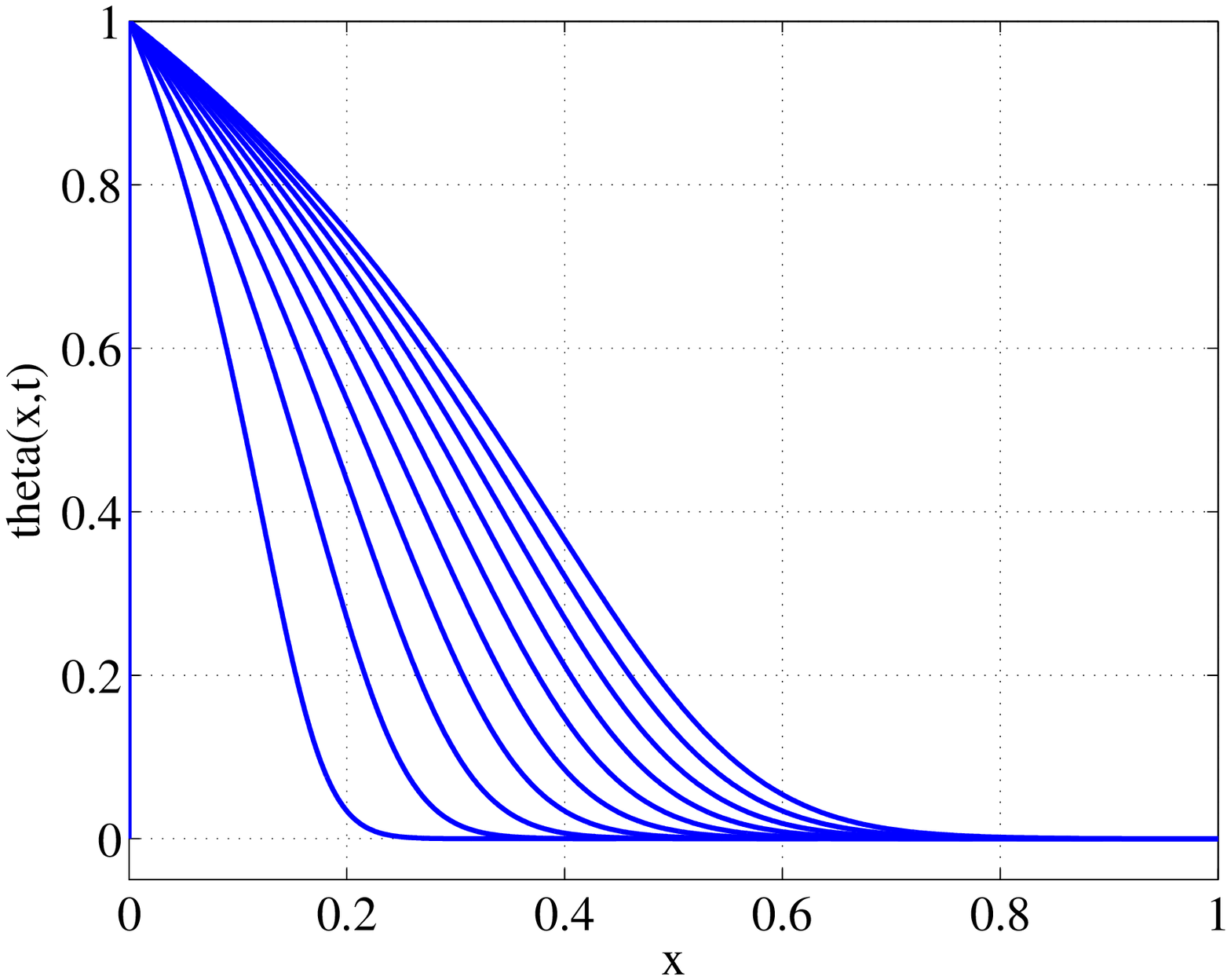}&
      \includegraphics[width=0.30\textwidth,clip]{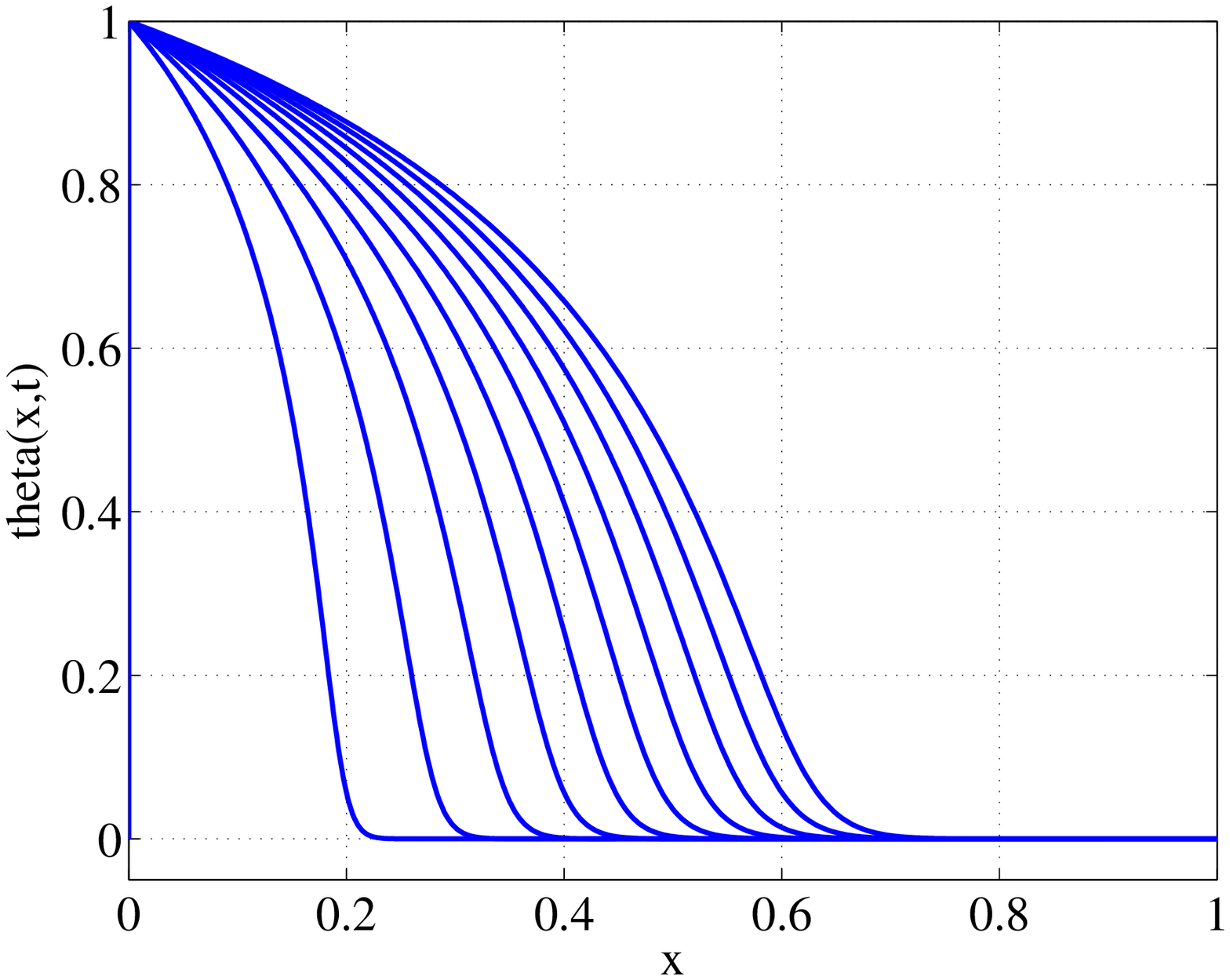}&
      \includegraphics[width=0.30\textwidth,clip]{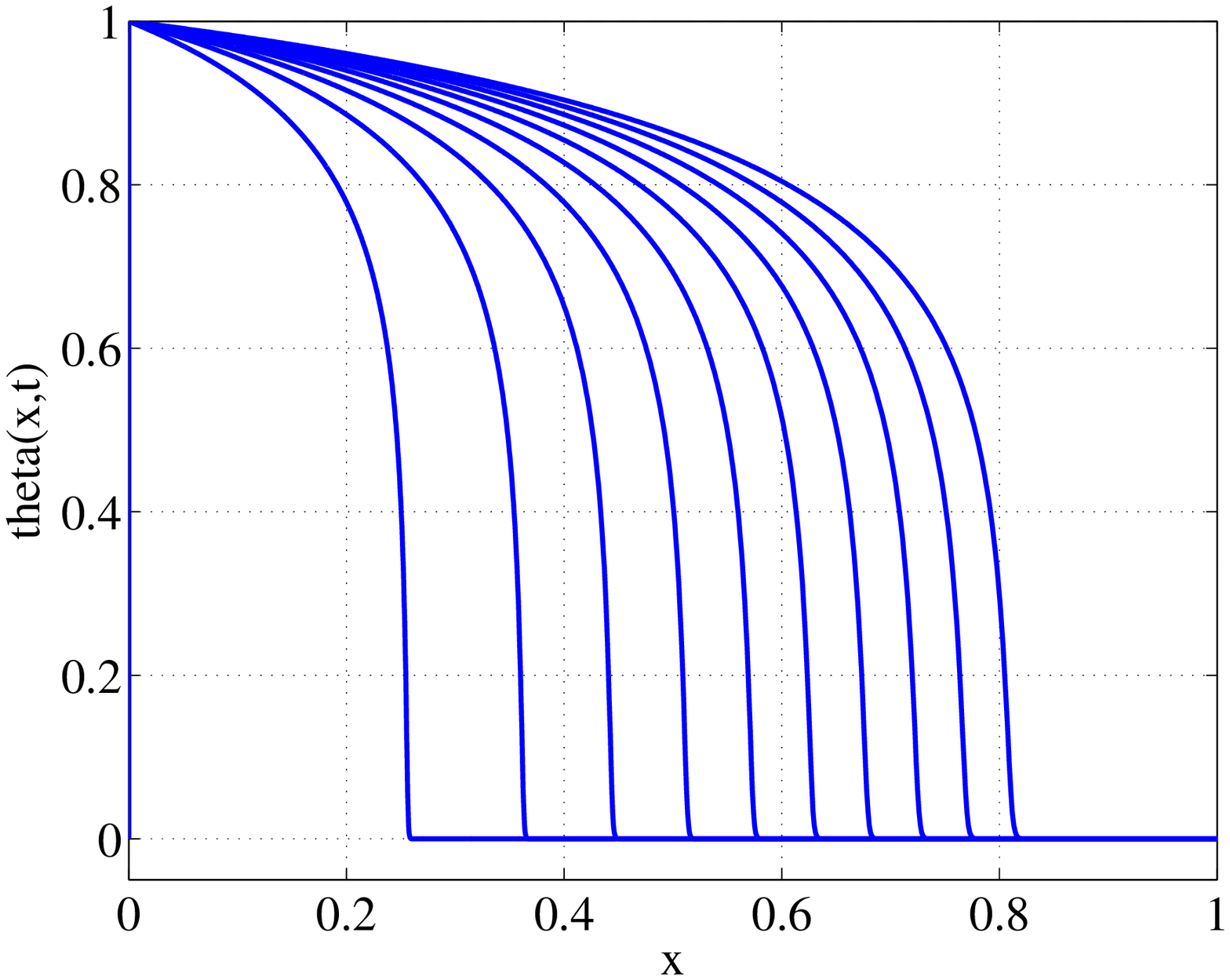}
    \end{tabular}
    \caption{Saturation for values of $\bbar = \beta\dtheta = 2$ (left), 4
      (middle) and 8 (right).  The top row shows the similarity variable
      $\Theta(y)$ obtained by solving {\protect Eqs.~\enprime{usim}}
      numerically.  The bottom row contains corresponding plots of the
      reduced saturation, $\sbar(x,t) = 1+(\log\Theta)/\bbar$, with
      time $t$ taken at ten equally-spaced points.  The inset at the top
      right is a zoomed-in view of the sharp corner.}
    \label{fig:gamma-all}
  \end{center}
\end{figure}

Before proceeding with the analysis, we present several plots that
illustrate the asymptotic behaviour of the self-similar solution for
large $\gamma$, computed via numerical simulations of the initial value
problem~\enprime{usim}.  We use a shooting algorithm wherein a value of
$\bbar$ is chosen, and then $\gamma=-\deriv{\Theta}{y}(0)$ is updated
iteratively until $\Theta$ is sufficiently close to the right hand
boundary value $\eps=e^{-\bbar}$.  This approach is similar to that
employed by others for the exponential diffusion
problem~\cite{crank-1975,yeh-franzini-1968}, and more details on our
shooting algorithm are given in Section~\ref{sec:numerical-algo}.  For
illustration purposes, we choose physical parameters corresponding to a
typical soil water uptake
experiment~\cite{reichardt-nielsen-biggar-1972}.  Taking limiting values
of $\theta_o=0.04$ and $\theta_i=0.43$, and an unsaturated diffusivity
${\beta}=20.5$, we obtain parameters \myrevision{$\dtheta=0.39$,
  $\bbar=8.0$ and $\eps =\exp(-8.0) = 3.4\times 10^{-4}$, which clearly
  satisfies the requirement $\eps\ll 1$}.  We remark that other
experiments on water transport in a wide range of porous media
(including soils, concrete, and other building materials) suggest that
allowable values of $\bbar$ are restricted to a fairly narrow interval
of $4\lessapprox\bbar\lessapprox 9$. The resulting numerical solution is
depicted in the two rightmost plots in Figure~\ref{fig:gamma-all}(c).
The upper plot displays the computed similarity solution $\Theta(y)$
while the lower plot shows the corresponding curves for reduced
saturation $\sbar(x,t)$ at ten equally-spaced time intervals, which are
determined from $\Theta(y)$ by transforming back to physical variables
Eq.~\en{boltzmann}.


The computed solution for $\bbar=8$ exhibits a well-defined wetting
front that manifests in the $\sbar$ plot as a narrow region with a
steep slope located immediately behind a sharp corner.  In the plot of
$\Theta$ on the other hand, the steep front has been eliminated by the
exponential stretching transformation \en{similar-sat}, but the wetting
front location is still identified with a sharp corner.  A zoomed-in
view of the wetting front is shown in the inset for $\bbar=8$ and
clearly indicates that the solution within the corner region transitions
rapidly to a small value of saturation, but this transition 
remains smooth.

Additional pairs of solution plots are given in
Figures~\ref{fig:gamma-all}(a,left) and \ref{fig:gamma-all}(b,middle)
for values of $\bbar=2$ and 4 respectively.  Clearly, taking $\bbar$
smaller (or equivalently $\eps$ larger) causes the wetting front to
exhibit a more gradual slope and a milder transition through the corner
region, and in the extreme case of $\bbar=2$ there is hardly any
evidence of a wetting front at all.  However, as we have indicated
above, most problems of physical interest correspond to values of $\bbar
\geqslant 4$ and this observation has an important effect on the
accuracy of our asymptotic solution derived in
Section~\ref{sec:asymptotic}.

\begin{figure}[tbhp]
  \psfrag{gamma}[bc][cc]{$\gamma$}
  \psfrag{beta}[cc][bc]{$\bbar$}
  \psfrag{ystar}[bc][cc]{$\ystar$}
  \begin{center}
    \leavevmode
    \begin{tabular}{cc}
      (a) Initial slope, $\gamma=-\deriv{\Theta}{y}(0)$ & 
      (b) Front location, $\ystar$ \\
      \includegraphics[width=0.385\textwidth,clip]{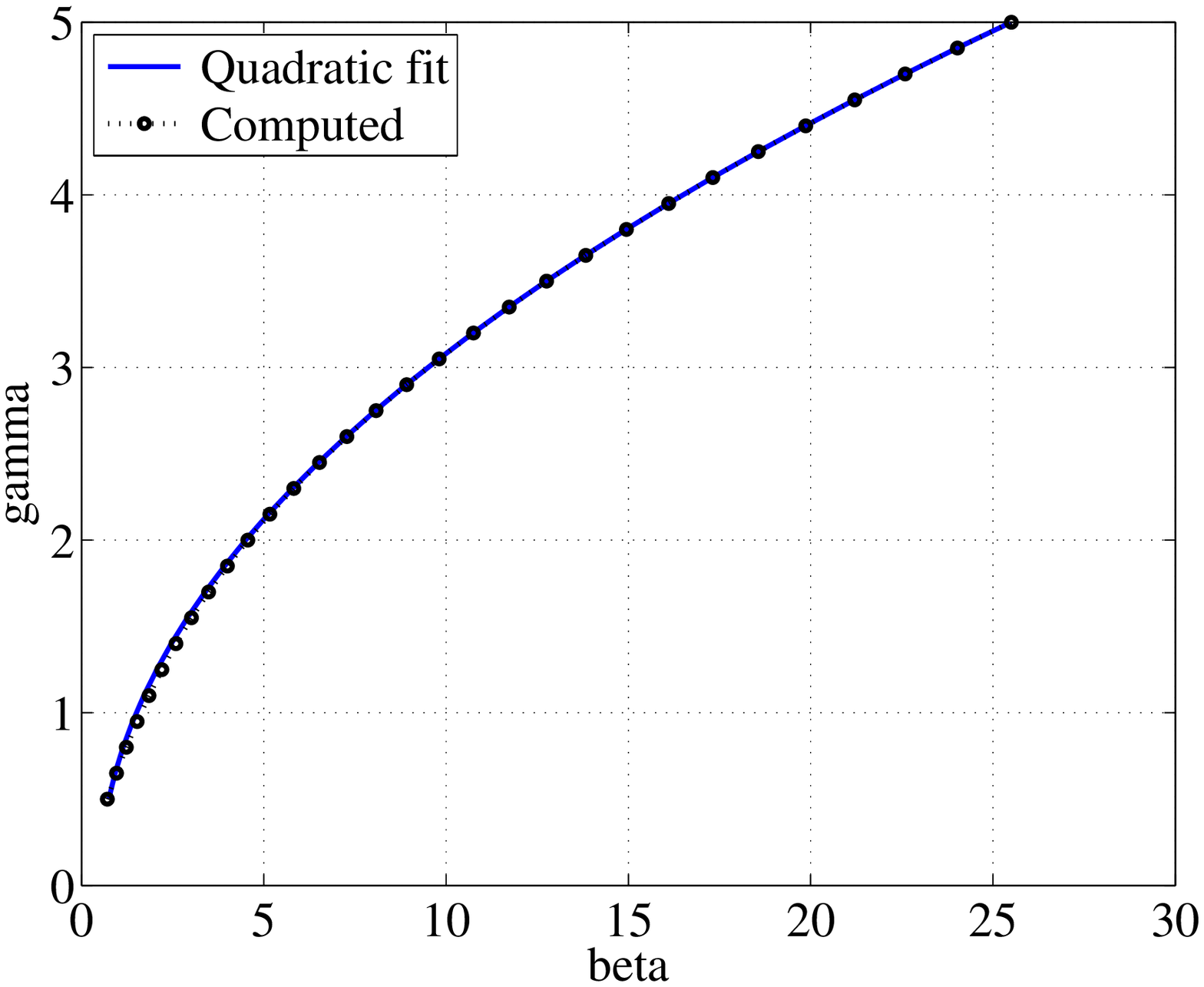} &
      \includegraphics[width=0.395\textwidth,clip]{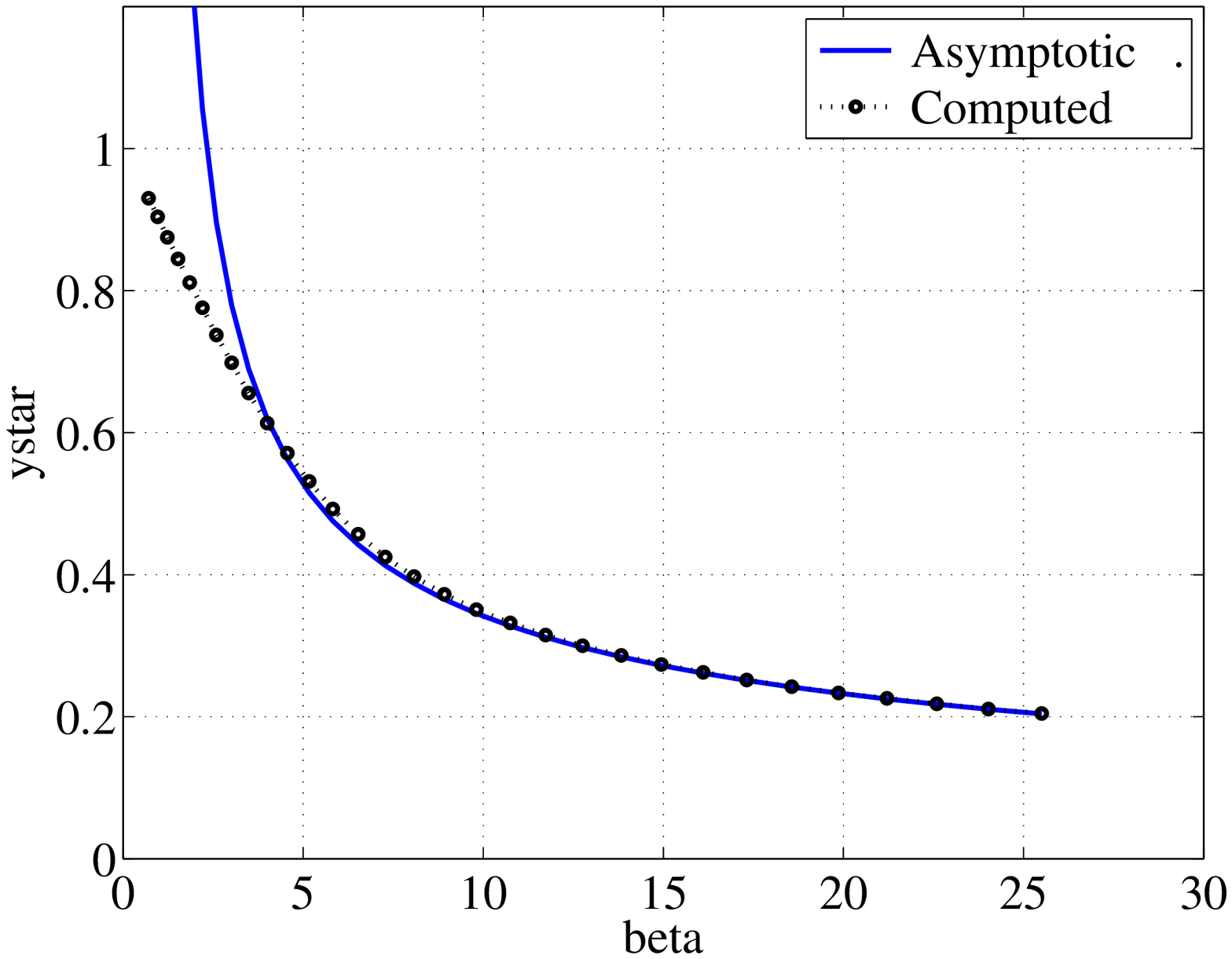} 
    \end{tabular}
    \caption{Left: Values of $\bbar$ and $\gamma=-\deriv{\Theta}{y}(0)$
      obtained from numerical simulations of {\protect
        Eq.~\enprime{usim}} with parameters chosen as in {\protect
        Figure~\ref{fig:gamma-all}}.  The quadratic fit
      $\bbar=\gamma^2+1/2$ is shown as a solid line for comparison
      purposes.  Right: The computed wetting front location (estimated
      using the point of maximum curvature) is depicted along with the
      two-term asymptotic approximation from {\protect
        Eq.~\en{ystar-asy}.}}
    \label{fig:betagam}
  \end{center}
\end{figure}

We next investigate in more detail the effect on the similarity solution
$\Theta(y)$ of changes in $\bbar$.  In particular,
Figure~\ref{fig:gamma-all} indicates that \myrevision{as $\bbar$
  increases the value of $\gamma$ likewise increases, which leads to a
  steeper initial slope and a corresponding shift} of the wetting front
location toward the origin along the $y$--axis.  A number of
additional simulations are performed for $\gamma$ lying in the interval
$[0.5, 5.0]$ and the corresponding values of $\bbar$ are plotted in
Figure~\ref{fig:betagam}(a).  By performing a least-squares polynomial
fit to the computed points, we obtain to a very good approximation the
quadratic polynomial fit
\begin{gather*}
  \bbar = - \log\eps \approx \gamma^2+\half, 
\end{gather*}
which when displayed in Figure~\ref{fig:betagam}(a) alongside the
computed data is nearly indistinguishable.  This relationship between
$\bbar$ and $\gamma$ will be verified later in
Section~\ref{sec:asymptotic} when it is derived as part of our
asymptotic solution.

It is evident from the plots in Figure~\ref{fig:gamma-all} that for
$\bbar$ sufficiently large there exists a sharp corner in the similarity
solution $\Theta(y)$ that can be identified with the wetting front
location in plots of saturation $\sbar$.  In contrast with the power-law
diffusion problem, where the wetting front is identified with a
discontinuity in the solution derivative, the exponential diffusion
problem exhibits a smooth transition through the front and so there is
no unique front position.  We therefore choose to identify the wetting
front location $\ystar$ with the point of maximum curvature in $\Theta$
which satisfies $\dderiv{\Theta}{y}(\ystar)=0$.  We provide a preview of
our asymptotic results in Figure~\ref{fig:betagam}(b), which depicts the
computed front location $\ystar$ as a function of $\bbar$, alongside our
our two-term asymptotic approximation of $\ystar$ derived later in
Section~\ref{sec:inner}.  Clearly, the analytical results are quite
accurate for $\bbar$ in the physical range.

\myrevision{\section{Multi-Layer Asymptotic Solution}}
\label{sec:asymptotic}

\subsection{Overview}

\begin{figure}[tbhp]
  \small
  \psfrag{u}{$\Theta$}
  \psfrag{=}{$\approx$}
  \centering
  {
    \psfrag{Th}{{\Large$\Theta$}}
    \psfrag{1}{{\Large 1}}
    \psfrag{0}{{\Large 0}}
    \psfrag{y}{{\Large $y$}}
    \psfrag{slope}{\textcolor{MyBlue}{\!\!\!\!$\deriv{\Theta}{y}(0)=-\gamma$}}
    \psfrag{inner}{{\bf \ Inner}}
    \psfrag{y=r/g}{$r = \gamma y$}
    \psfrag{Th=u(r)}[B][b]{$\Theta=\vinn(r)$}
    \psfrag{mid-range}{{\bf \hspace*{0.2cm}\myrevision{Intermediate}}}
    \psfrag{y=1/g+s/g^3}{$s=\gamma^2(\gamma y-1)$}
    \psfrag{Th=v(s)/g^4}{\hspace*{0.3cm}${\displaystyle \Theta=\frac{v(s)}{\gamma^2}}$}
    \psfrag{(s<1/2)}[b][t]{{\footnotesize \hspace*{0.3cm}$\left(s<\half\right)$}}
    \psfrag{Th=w(s)*Th8}{{$\Theta=\eps w(s)$}}
    \psfrag{(s>1/2)}[b][t]{{\footnotesize \hspace*{0.2cm}$\left(s>\half\right)$}}
    \psfrag{y=1/g}{\textcolor{MyBlue}{${\displaystyle y=\frac{1}{\gamma}}$}}
    \psfrag{(r=1)}[t][t]{\ \textcolor{MyBlue}{{\footnotesize $(r=1)$}}}
    \psfrag{y*=1/g+1/2g^3}{\hspace*{-0.3cm}\textcolor{MyRed}{${\displaystyle
          \ystar\approx \frac{1}{\gamma} + \frac{1}{2\gamma^3} +
          \frac{11}{12\gamma^5}}$}}  
    \psfrag{(s=1/2)}[t][t]{\hspace*{0.7cm}\textcolor{MyRed}{\footnotesize
        $\left(\sstar\approx \half+\frac{11}{12\gamma^2}\right)$}}
    \psfrag{outer}{{\bf \ Outer}}
    \psfrag{Th=Th8(1+g(y))}{$\Theta=\eps(1+g(y))$}
    \psfrag{Th8}{$\eps = e^{-\bbar}$}
    \includegraphics[width=0.55\textwidth]{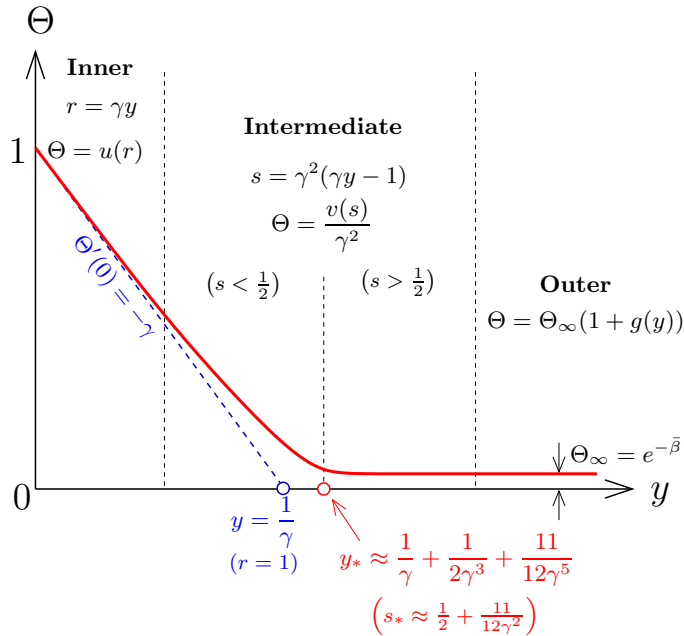}
  }  
  \caption{The similarity solution $\Theta(y)$ is separated into several
    regions: an inner region on the left consisting of a nearly linear
    solution $\vinn(r)$; an outer region on the right where the solution
    is nearly constant; and an intermediate region near the wetting
    front $y=\ystar$ within which an \myrevision{intermediate-range}
    solution of locally high curvature connects the inner and outer
    solutions.}
  \label{fig:multilayer}
\end{figure}

In this section, we will derive the asymptotic form of the solution for
large $\gamma$. The preliminary numerical results already shown in
Figure~\ref{fig:gamma-all} suggest that the solution for large $\bbar$
(and hence large $\gamma$) can be separated into four regions or
layers:
\begin{enumerate}
  \renewcommand{\labelenumi}{(\alph{enumi})}
\item An linear \emph{inner solution} in which $\Theta(y)$ is close to
  linear and has gradient close to $1/\gamma$.  This region extends over
  a range of $y$ values lying between 0 and slightly below $1/\gamma$
  where the solution is determined by the initial conditions at $y=0$.
\item A nearly constant \emph{outer solution} on the right of $\ystar$,
  where the solution is determined by the far field condition as $y \to
  \infty$. Here the solution approaches $\eps$ for $y$ sufficiently
  large (in fact, for $y$ only a little greater than $\ystar$).
\item An \emph{intermediate-range solution} containing the wetting front
  at $\ystar$ and of width $\Delta y = {\cal O}(1/\gamma^3)$.  In this
  region we see the rapid transition through the corner of the wetting
  front and the solution has (locally) and exponentially (in $\gamma$)
  large curvature. In the intermediate range, we see a transition from
  the linear solution, which is of order $1/\gamma^2$ when $y =
  1/\gamma$, to one exponentially small in $\gamma$ for $y > \ystar =
  1/\gamma + 1/2\gamma^3 + \order{1/\gamma^5}$.  It will prove
  convenient to divide this intermediate layer into a \emph{left-range}
  and a \emph{right-range} to capture this behaviour to high order.
\end{enumerate}
These regions are depicted schematically in Figure~\ref{fig:multilayer}
in terms of the similarity variable $\Theta(y)$.  We expect that our
asymptotic approximation will be inaccurate for small values of $\bbar$
when the inner solution is more curved behind the front, but will
improve as $\bbar$ increases and the inner solution becomes closer to
linear.
%
%

Before going into the details, we first summarise the main result for
the wetting front location
\begin{gather}
  \ystar = \frac{1}{\gamma} + \frac{1}{2\gamma^3} +
  \frac{11}{12 \gamma^5} + \order{\frac{1}{\gamma^7}}. 
  \label{eq:ystar-prelim}
\end{gather}
By taking $\ystar$ as the point of maximum curvature in saturation
$\Theta(y)$, we may then approximate the values of both saturation
and curvature at the wetting front by
\begin{gather}
  \Theta(\ystar) \approx e\eps \quad 
  \text{and}\quad
  \dderiv{\Theta}{y}(\ystar) \approx \frac{1}{\gamma^2}
  e^{\gamma^2-1/2}. 
\end{gather}
We will show further that 
\begin{gather}
  \Theta\left(\frac{1}{\gamma}\right) = \frac{1}{2\gamma^2} +
  \frac{b-\log(\gamma)}{\gamma^4} + \order{\frac{1}{\gamma^6}}
  \quad \text{where} \quad
  b = \frac{11}{12} - \half\log{2}. 
\end{gather}
Finally, the physical parameter $\bbar$ can be expressed as 
\begin{gather}
  \bbar = -\log\eps = \gamma^2 + \half + \frac{\alpha_3}{\gamma^2} +
  \order{\frac{1}{\gamma^4}}, 
\end{gather}
where $\alpha_3$ is a constant that we can estimate numerically as
$\alpha_3 = 1/12$. It follows immediately that for large $\bbar$ we
have
\begin{gather}
  \gamma = \bbar^{1/2} - \frac{1}{4} \bbar^{-1/2} -
  \left(\frac{\alpha_3}{2} - \frac{1}{32}\right) \bbar^{-3/2} +
  \order{\bbar^{-5/2}}.
\end{gather}

\subsection{Inner Solution}
\label{sec:inner}

\begin{figure}[tbhp]
  \centering
  \psfrag{y}[cc][bc]{$y$}
  \psfrag{Theta\(y\)}[bc][cc]{$\Theta(y)$}
  \psfrag{       Inner \(1 term\)}{{\scriptsize \!\!\!\!\!\!\!\! $\vinnz$}}
  \psfrag{       Inner \(2 terms\)}{{\scriptsize \!\!\!\!\!\!\!\! $\vinnz+\vinno/\gamma^2$}}
  \psfrag{       Exact}{{\scriptsize \!\!\!\!\!\!\!\! Exact}}
  \begin{tabular}{ccc}
    (a)\ \ $\bbar=2$,\ \ $\gamma=1.166$ 
    & 
    (b)\ \ $\bbar=4$,\ \ $\gamma=1.850$ 
    & 
    (c)\ \ $\bbar=8$,\ \ $\gamma=2.736$ 
    \\
    \includegraphics[width=0.30\textwidth,clip]{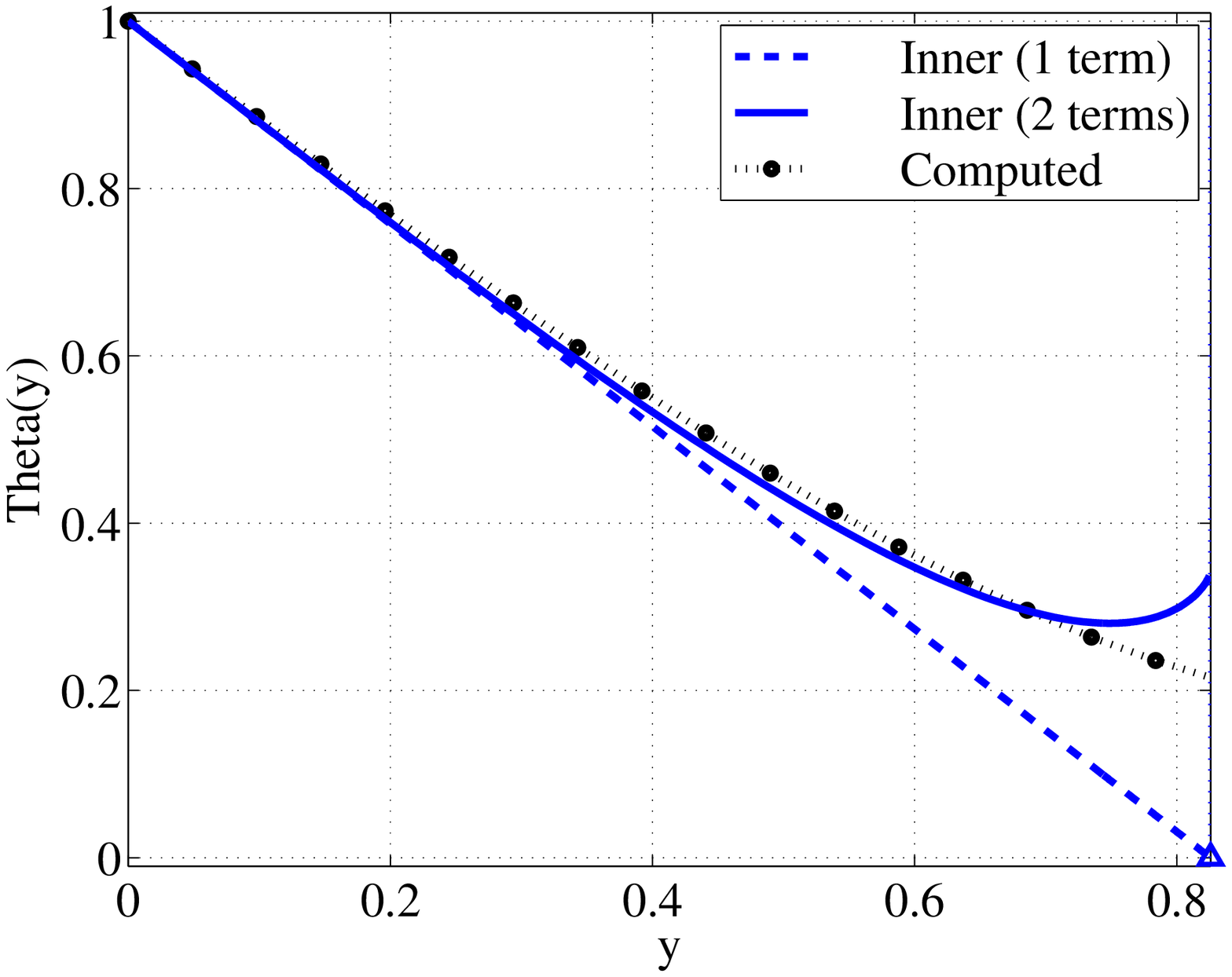} &
    \includegraphics[width=0.30\textwidth,clip]{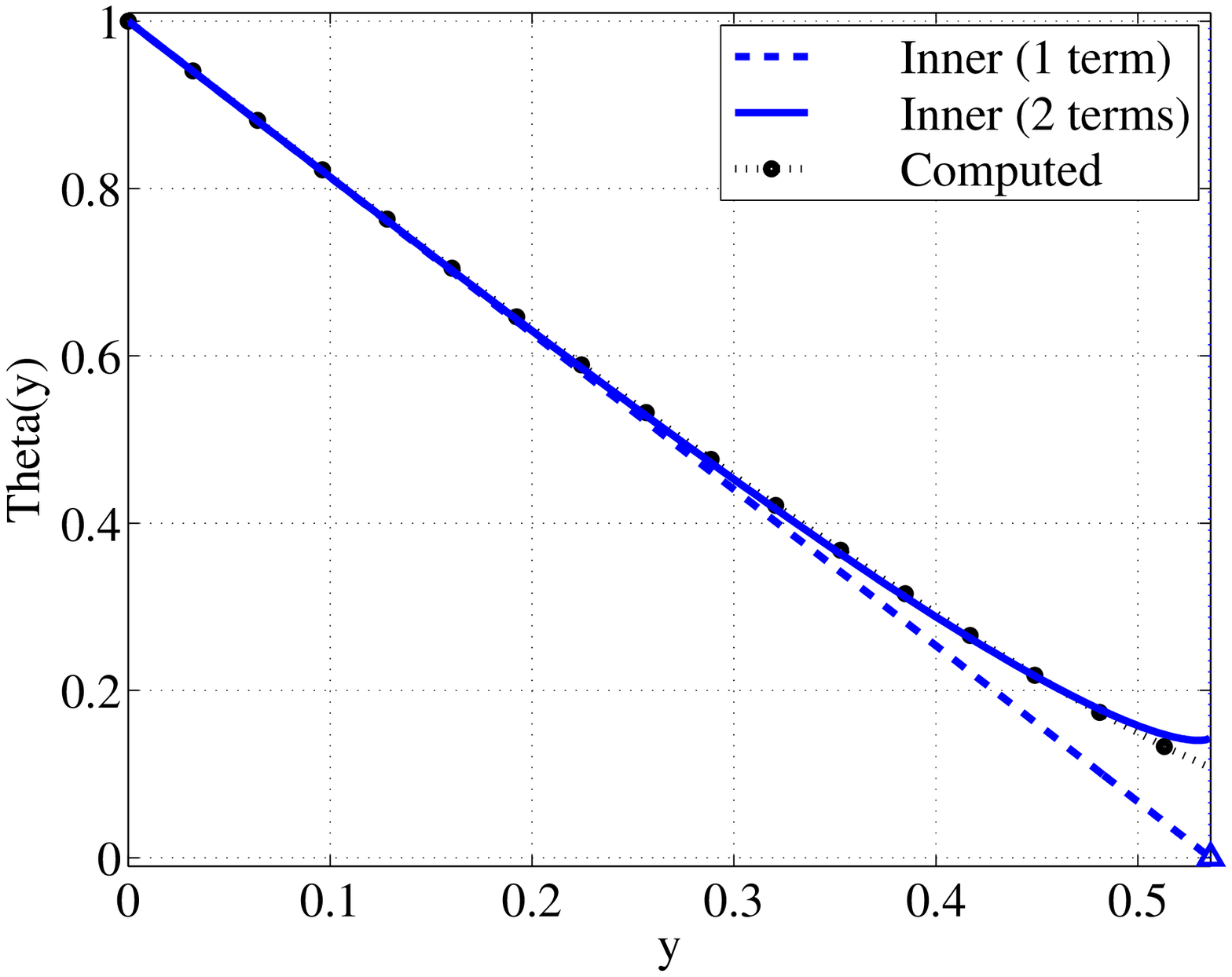} &
    \includegraphics[width=0.30\textwidth,clip]{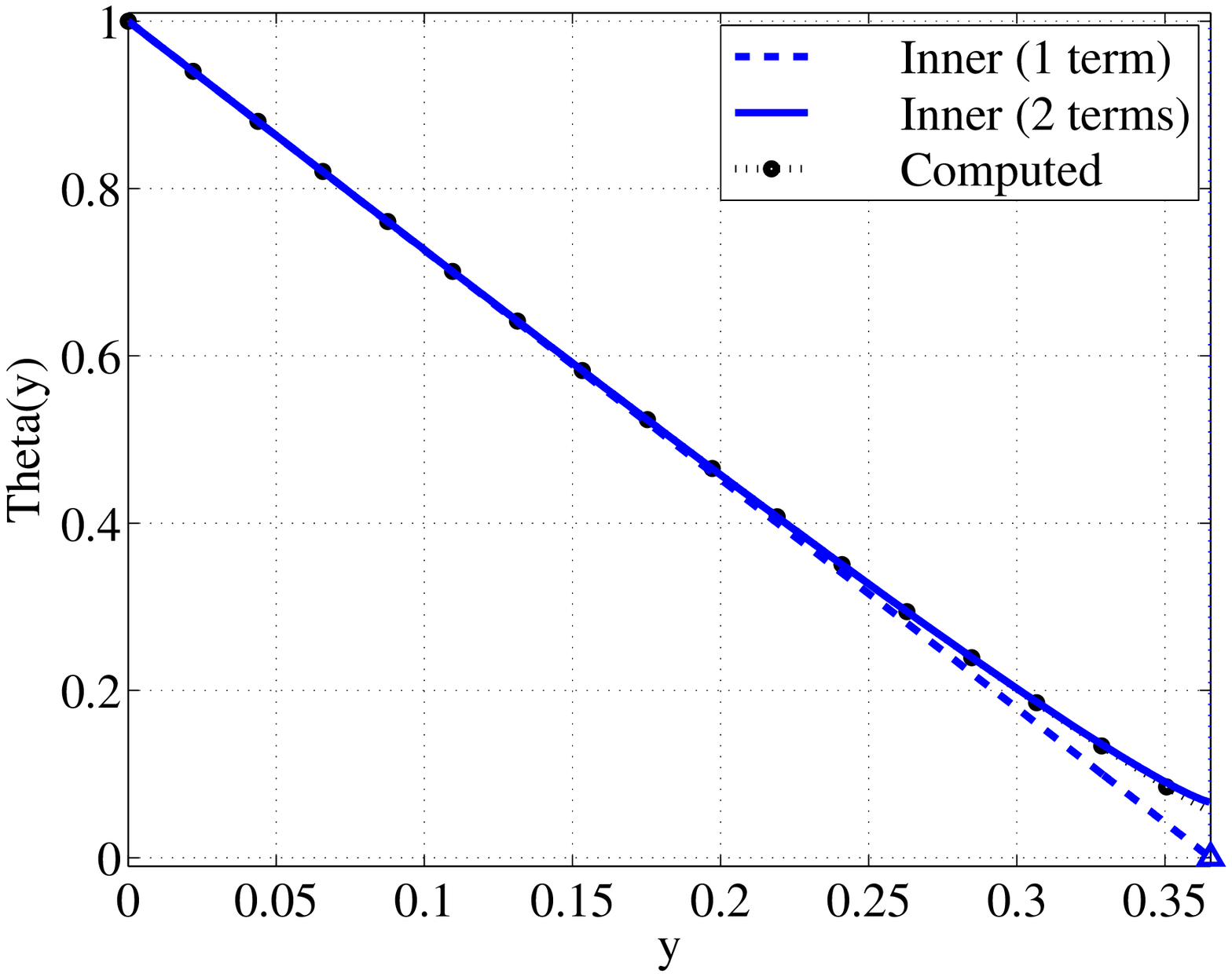}
  \end{tabular}
  \caption{The inner solution $\vinn(r)$ for values of $\bbar=2$, 4 and
    8, showing the leading order approximation $\vinnz$, first order
    correction $\vinnz+\vinno/\gamma^2$, and numerical solution of the
    ODE initial value problem.  The point where the approximation
    $\vinnz$ touches the $y$--axis (\textcolor{MyBlue}{\small
      $\triangle$}) corresponds to the leading order estimate of the
    wetting front location, $\ystar\approx 1/\gamma$.}
  \label{fig:vinn}
\end{figure}

On the interval $0 \leqslant y < \ystar$, we know that the solution
slope is initially $-\gamma$ and so it is natural to introduce a new
independent variable
\begin{gather*}
  r = \gamma y,
\end{gather*}
and define the inner solution as $\vinn(r) = \Theta(y)$.  Rescaling
Eqs.~\enprime{usim} yields
\begin{subequations}\label{eq:inner}
  \begin{align}
    \vinn \,\dderiv{\vinn}{r} &=
    -\frac{r\,\deriv{\vinn}{r}}{\gamma^2},\label{eq:inner-v}\\  
    \vinn(0) &= 1,\label{eq:inner-bc1}\\
    \deriv{\vinn}{r}(0) &= -1\label{eq:inner-bc2}.
  \end{align}
\end{subequations}%
%
%
The $1/\gamma^2$ factor on the right hand side of \en{inner-v} suggests
that if $\gamma$ is large then we should use an asymptotic expansion of
the form
\begin{subequations}\label{eq:vinn}
  \begin{gather}
    \vinn(r) = \vinnz(r) + \frac{\vinno(r)}{\gamma^2} + 
    \frac{\vinnt(r)}{\gamma^4} + 
    \order{\frac{1}{\gamma^6}}.
    \label{eq:vinn-series}
  \end{gather}
  After substituting this expression into the ODE and boundary
  conditions for $u$, we obtain the following sequence of initial value
  problems up to $\order{\gamma^{-4}}$:
  \begin{alignat*}{8}
    \vinnz\,\dderiv{\vinnz}{r} & = 0, 
    & \qquad & 
    \vinnz(0)\, & = 1, & \quad & 
    \deriv{\vinnz}{r}(0) & = -1,\\
    \vinnz\,\dderiv{\vinno}{r} & = -r\deriv{\vinnz}{r}, 
    & &
    \vinno(0) \, & = 0, & \quad & 
    \deriv{\vinno}{r}(0) & = 0,\\
    \vinnz\,\dderiv{\vinnt}{r} & = -\vinno\dderiv{\vinno}{r} -
    r\deriv{\vinno}{r}, 
    & & \vinnt(0) \, & = 0, & & 
    \deriv{\vinnt}{r}(0) & = 0.
  \end{alignat*}
  These problems can be integrated successively to obtain
  \begin{align}
    \vinnz & = 1-r,
    \label{eq:vinner0} \\
    \vinno & = \half  - \half (1-r)^2 + (1-r)\, \log(1-r),  
    \label{eq:vinner1}\\
    \vinnt & = 
    \frac{17}{12} - \frac{3}{4}(1-r) - \frac{3}{4}(1-r)^2
    + \frac{1}{12}(1-r)^3 + \left(2-\frac{3}{2}r\right) \,
    \log(1-r),  
    \label{eq:vinner2}
  \end{align}
\end{subequations}
which after substitution into \en{vinn-series} gives the inner solution
to $\order{\gamma^{-4}}$.  The $\order{1}$ and $\order{\gamma^{-2}}$
solutions are depicted in Figure~\ref{fig:vinn} for values of $\bbar=2$, 4
and 8.  Plots of the numerical solution of Eqs.~\enprime{usim}, computed
using the shooting algorithm described in Section~\ref{sec:numerical}
(and labeled ``Exact'') have also been included for comparison
purposes.  

We now investigate more carefully the validity of these asymptotic
expressions. It is immediately clear that for asymptotic regularity we
should have
\begin{gather*}
  \vinnz(r) \gg \frac{\vinno(r)}{\gamma^2}.
\end{gather*}
Taking the leading order term on each side of this equation yields
\begin{gather}
  (1 - r) \gg \frac{1}{2 \gamma^{2}}.
\end{gather}

Finally, we discuss how the first two terms in the asymptotic expansion
\en{vinn} may be used to give a first estimate the wetting front
location $\ystar$. To this end, we set $\vinnz+\vinno/\gamma^2=0$ and
neglect terms that are quadratic and logarithmic in $(1-r)$ to obtain
\begin{gather}
  \rstar \approx 1 + \frac{1}{2 \gamma^2} \qquad
  \mbox{or} \qquad \ystar \approx \frac{1}{\gamma} +
  \frac{1}{2\gamma^3},  
  \label{eq:ystar-asy}
\end{gather}
which yields the first two terms in the front location
\en{ystar-prelim}. We also note that setting $r=1$ yields the leading
order estimates
\begin{gather}
  u(1) \approx \frac{1}{2 \gamma^2} \qquad \text{and} 
  \qquad u^{\prime}(1) = -1,
  \label{eq:u0est}
\end{gather}
which will be useful in later scaling arguments.  One of the primary
aims of the more careful asymptotic matching performed in
Section~\ref{sec:mid-mid} is to derive a correction to
Eq.~\en{ystar-asy} that gives more refined estimates of $u(1)$ and 
$\rstar$ to verify the above rough calculation.

\subsection{Intermediate Layer}
\label{sec:midrange}

\subsubsection{Preliminary Estimate of Leading Order Asymptotics and Solution Scalings} 
\label{sec:prelim}

As $r \to 1$ we start to see the rapid transition between the linear
solution and the exponentially small solution.  We note that a standard
application of the maximum principle in the intermediate region shows
that $u(r) > 0$ and $\deriv{u}{r} < 0$ for all $r$, so that $u \to
\uinf \equiv \eps$ and $\deriv{u}{r} \to 0$ as $r \to \infty$.
This transition occurs over a narrow layer, which we will see is of
width $\order{1/\gamma^2}$.  Before getting into the detailed analysis
of intermediate solution, we begin with some simple calculations that
aim to establish the leading order asymptotic form of the solution in
intermediate layer close to the wetting front at $\rstar$ and also to
estimate the width of this layer. These calculations will guide our
more sophisticated calculations given later.  Using the scalings in the
inner region, we have
\begin{gather}
  u\dderiv{u}{r} = -\frac{r \deriv{u}{r}}{\gamma^2}, 
  \quad u(0) = 1,
  \quad \deriv{u}{u}(1^-) = -1, 
  \quad u(\infty) = \uinf \equiv \eps.
\end{gather}
If we now divide both sides of the ODE by $r u$, then integrate and
apply the boundary conditions, we obtain
\begin{gather}
  -\log(\uinf) = \gamma^2 \int_{0}^{\infty} \frac{\dderiv{u}{r}(r)}{r} \,
  dr. \label{eq:intest}  
\end{gather}
In the case when $\gamma$ is large, our previous estimates of the inner
solution, combined with the numerical calculations, imply that:
\begin{center}
  (a) $\deriv{u}{r} \approx -1$ if $0 \leqslant r < \rstar$;\quad
  (b) $\deriv{u}{r} \approx 0$ if $r > \rstar $; \quad
  (c) $\dderiv{u}{r} \approx 0$ if $r \ne \rstar$; \ \ and\ \ 
  (d) $\dderiv{u}{r}(\rstar) \gg 1$.
\end{center}
It follows from these observations that
\begin{gather*}
  \gamma^2 \int_{0}^{\infty} \frac{\dderiv{u}{r}}{r} \, dr
  \approx 
  \gamma^2\int^{\rstar^+}_{\rstar^-} \frac{\dderiv{u}{r}}{r} \, dr
  \approx
  \gamma^2 \Big[ \frac{\deriv{u}{r}}{r} \Big]_{\rstar^-}^{\rstar^+} 
  \approx \frac{\gamma^{2}}{\rstar},
\end{gather*}
which can be combined with \en{intest} and the inner solution estimate
$\rstar = 1 + {\cal O}(1/\gamma^2)$ to obtain to leading order
\begin{gather}
  -\log \uinf \approx \gamma^2.
\end{gather}
While this estimate for $\uinf$ is not precise, we can use this
simple calculation as a guide in the subsequent development of the
asymptotic theory.

For our next estimate, we make the assumption that as $u^{\prime\prime}$
is exponentially large only close to $\rstar$, the leading order
behaviour of the integral in \en{intest}) can be approximated by
freezing the value of $r=r*$ in the denominator of this integral.  By
integrating and applying the previous result, we obtain the leading
order equation
\begin{gather}
  -\frac{\gamma^{2}  u'}{\rstar} = \log \left( \frac{u}{\uinf} \right).
  \label{eq:lead1}
\end{gather}
This equation can be integrated exactly when $r > 1$ to obtain
\begin{gather}
  \rstar \left( r -1 \right)  = \gamma^{2} \uinf \left[ \li\left(
      \frac{u(1)}{\uinf} \right) - 
    \li\left( \frac{u}{\uinf} \right) \right],
  \label{eq:midmatch1}
\end{gather}
where
\begin{gather*}
  \li(z) = \int_0^z \frac{dt}{\log t}
\end{gather*} 
is the \emph{logarithmic integral function}.  In order to determine the
leading order behaviour in the intermediate layer $r \approx 1$ and to
match with the inner and outer solutions, we consider the limit
$u(r)/\uinf \to 1$.  From \en{u0est}, we have that $u(1)\approx
1/2\gamma^2$ and $\uinf \approx \exp(-\gamma^2)$, so that
$u(1)/\uinf \gg 1$ for large $\gamma$.  To calculate the solution
thus requires estimates of the logarithmic integral $\li(z)$ in the two
limits $z\to 1$ and $z\to\infty$, which are well-known to be
\begin{alignat}{3}
  \li(z) &= \eulergam + \log\log z + \order{\log z} 
  \quad &&\text{as} \quad  z\to 1, \label{eq:li-asy1} \\
  \li(z) &= \frac{z}{\log z} + \frac{z}{\log^2 z} + \order{\frac{z}{\log^3 z}} 
  \quad &&\text{as} \quad  z\to \infty, \label{eq:li-asy2}  
\end{alignat}
where $\eulergam=0.5772156649\dots$ is the Euler-Mascheroni constant.\ \ %
%
%
Using the expansion \en{li-asy2} for large $z$, we obtain 
\begin{align}
  \gamma^2 \uinf \li\left( \frac{u(1)}{\uinf} \right) 
  & = \gamma^2 \frac{u(1)}{\log(u(1)) - \log(\uinf)} + 
  \order{\frac{\gamma^2 u(1)}{\log^2(\uinf)}} \notag\\ 
  & =  \frac{1}{2 \gamma^2 } + \order{\frac{1}{\gamma^4}} \notag\\ 
\end{align}

If we now let $u/\uinf$ be close to one, so that we are very close
to the wetting front, and then combine the various expressions above, we have
\begin{gather}
  \rstar ( r - 1 ) = \frac{1}{2 \gamma^2} +  \order{\frac{1}{\gamma^4}}
  - \gamma^2 \uinf \Big( 
    \eulergam + \log \log (u/\uinf) \Big).
\end{gather}
Because $\rstar = 1 + \order{1/\gamma^2}$, this is equivalent to
\begin{gather}
  u = \uinf \exp \exp \left( -\gamma^{-2} \uinf^{-1} \left( r
      - 1 - \frac{1}{2 \gamma^2} - \order{\frac{1}{\gamma^4}}  
    \right) - \eulergam \right).
  \label{eq:inter1}
\end{gather}
If (as consistent with our earlier estimates) we identify the term $1 +
\frac{1}{2 \gamma^2} + \order{\frac{1}{\gamma^4}}$ in this expression
with the wetting front at $\rstar$, we then have
\begin{gather}
  u = \uinf \exp \exp \Big( -\frac{r-\rstar}{\gamma^2 \uinf}
  - \eulergam \Big). 
  \label{eq:inter2}
\end{gather}
The structure of the intermediate layer is now clear.  Because
$\uinf^{-1} \approx e^{\gamma^2}$ is very large indeed, there is a
very rapid transition from the inner solution where $u^\prime \approx
-1$ to the solution $u = \uinf, u^\prime=0$ as we pass through the
wetting front, over a range $\Delta r \approx \gamma^2 e^{-\gamma^2}$.

\subsubsection{Left Intermediate Solution}

Guided by the above, we now refine the asymptotic calculation in the
intermediate layer, in particular for the range $r < 1 + 1/2\gamma^2 +
{\cal O}(1/\gamma^4)$, to get a more precise estimate of the location of
the wetting front and of the value of $\uinf$. This calculation is
rather technical, but allows us to obtain higher order asymptotic
estimates which can then be compared with the numerical
calculations. Within this range, we first determine an appropriate
rescaling of the spatial variable $r$ and dependent variable $u$ by
studying the form of the inner solution as $r\to 1^-$.  In this limit,
the dominant terms in the inner expansion \en{vinn} are
\begin{gather}
  \vinn(r) = 1 - r + \frac{1}{\gamma^2} \left( \half + \dots
  \right), 
  \label{eq:vmid-leading}
\end{gather}
Consequently, we make introduce a new independent variable $s$
that is $\order{1}$ and is given by 
\begin{gather}
  r = 1 + \frac{s}{\gamma^2} 
  \label{eq:ymid}
\end{gather}
We note that taking $s$ large and negative (for example ${\cal
  O}(\gamma)$) allows us to match to the inner solution.  Furthermore,
the results of the previous section show that there is very rapid
convergence to the outer solution, with $u \to \uinf$, for $s
> 1/2 + {\cal O}(1/\gamma^2).$ We are motivated by considerations of
matching to define a new \myrevision{intermediate-range} saturation
variable $v(s)$ via
\begin{gather*}
  u(r) = \frac{v(s)}{\gamma^{2}}, 
\end{gather*}
so that
\begin{gather*}
  u^\prime(1^-) = v^\prime(-\infty) \quad \mbox{and} \quad  u(1^-) \approx
  \frac{1}{2 \gamma^2}  \quad \mbox{imply} \quad v(0) \approx
  \frac{1}{2}.
\end{gather*}
We then have \myrevision{to leading order}
\begin{gather}
  v \, v^{\prime\prime} = - \frac{\deriv{v}{s}}{\gamma^{2}}\, 
  \left( 1 + \frac{s}{\gamma^2} \right) .
  \label{eq:vmid-ode}
\end{gather}%
%
%
As $s$ increases the \myrevision{intermediate-range} solution will match
rapidly to the outer solution so that
\begin{gather}
  v(s) \to \vinf \defeq \gamma^2  \uinf, \quad v^\prime(s) \to 0 \quad
  \mbox{if} \quad s >  s_\ast = 1/2 + {\cal O}(1/\gamma^2) .
\end{gather}
%
%
Integrating Eq.~\en{vmid-ode} over the interval $[s,s_\ast^+)$ yields
\begin{gather*}
 \gamma^{-2} \Big( \log(v(s)) - \log(\vinf) \Big) = \int_{s}^{s_\ast}
  \frac{\dderiv{v}{s} \, ds}{(1 + s/\gamma^{2})}.
\end{gather*}
We will assume that the the range of $s$ is restricted to 
\begin{gather*}
  s/\gamma^2 \ll 1, 
\end{gather*}
so that the denominator in the integrand may be expanded as a geometric
series and
\begin{gather*}
 \gamma^{-2} \Big( \log(v(s)) - \log(\vinf) \Big) 
  =  \int_{s}^{s_\ast} \dderiv{v}{s}
  \left( 1 - \frac{s}{\gamma^2} + \order{\frac{s^2}{\gamma^4}}
  \right) \, ds.
\end{gather*}
Using integration by parts and applying the far field condition
$\deriv{v}{s}\to 0$ for $s_\ast^+$, we find that
\begin{gather*}
  \gamma^{-2} \Big( \log(v) - \log(\vinf) \Big)
  = - \deriv{v}{s} \left( 1-\frac{s}{\gamma^2} \right)
  + \frac{\left(\vinf - v \right)}{\gamma^2}
  + \order{\frac{1}{\gamma^4}},  
\end{gather*}
which can be rearranged to obtain
\begin{gather}
  \deriv{v}{s} = \gamma^{-2} \Big( \log(\vinf) - \log(v) \Big)
  + \frac{s\deriv{v}{s}}{\gamma^2} 
  + \frac{1}{\gamma^2} \left(\vinf - v\right)  
  + \order{\frac{1}{\gamma^4}}.  
  \label{eq:dvmid-approx}
\end{gather}

We now carefully consider the solution of this equation for $s< 1/2$ and
match it to the inner solution.  This calculation is rather technical
and involved, however it furnishes very precise information about the
value of $\eps$ and the wetting front location. To do this we develop
an asymptotic expression of the form
\begin{gather}
  v(s) = v_0(s) + \frac{v_1(s)}{\gamma^2}  +
  \frac{v_2(s)}{\gamma^4} + \ldots  \quad \text{when} \quad s < \half +
  {\cal O}(1/\gamma^2). 
  \label{eq:vmid-asy}
\end{gather}
A subtle feature of this expansion follows from the coexistence of
the logarithmic and polynomial expressions in both the intermediate
layer expansion and the inner solution. This means that in order
to capture the solution behaviour it is necessary that $v_1(s)$ has
terms of both $\order{1}$ and $\order{\log\gamma}$. However, provided
that $\gamma$ is large, this means that we still have $v_0 \gg
v_1/\gamma^2 \gg v_2/\gamma^4$ so that the formal nature of the
asymptotic series is preserved.  We assume further, motivated by the
previous analysis, that
\begin{align}
  v(0) &= \frac{1}{2} + a \frac{\log\gamma}{\gamma^2} +
  \frac{b}{\gamma^2} + \order{\frac{1}{\gamma^4}},  
  \label{eq:v0-mid}
  \\
  \log\vinf &= -\gamma^2 + c\log\gamma + d + \order{\frac{1}{\gamma^2}}, 
  \label{eq:vinf-mid}
\end{align} 
and we will determine explicit values for the constants $a$, $b$, $c$,
$d$.  We next substitute the expressions \en{vmid-asy}--\en{vinf-mid}
into \en{dvmid-approx} and solve the equations arising at each order.

Considering first the leading terms of $\order{1}$, we find that
\begin{gather*}
  \deriv{v_0}{s} = -1, \qquad 
  v_0(0) = \half,
\end{gather*}
which has solution
\begin{gather}
  v_0(s) = \half - s.
  \label{eq:vmid0-soln}
\end{gather}
Next, take the terms of $\order{1/\gamma^2, \log\gamma/\gamma^2}$ for
which we first need to expand $\log v$ from Eq.~\en{dvmid-approx} as a
series in $\gamma$
\begin{gather*}
  \log v = \log \left(  \left( \half - s \right) + v_1/\gamma^2 +
    \cdots \right) = \log\left( \half - s \right) +
  \order{\frac{1}{\gamma^2}},
\end{gather*}
where we have used the fact that $v_0>0$ when $s<1/2$.  The
$\order{1/\gamma^2,\log\gamma/\gamma^2}$ terms in the intermediate layer
expression then become
\begin{gather*}
  \deriv{v_1}{s} = c \; \log\gamma + \left( d-\half \right) - \log
  \left( \half - s \right), 
\end{gather*}
which can be integrated to obtain
\begin{gather}
  v_1(s) = c \; s \log\gamma + \left( d-\half \right) s 
  + a \log\gamma + b + \left( \half-s \right) \log \left( \half -
    s \right) + s - \half \log 2.
  \label{eq:vmid1}
\end{gather}

We now match the intermediate range expansion $v= v_0 + v_1/\gamma^2 +
\order{1/\gamma^4}$ to the inner solution when $s < 0$ and $|s| =
\order{\gamma} \gg 1$. Making use of the expansion
\begin{gather*}
  \left(\half-s\right) \log\left(\half-s\right) = 
  \left(\half - s\right) \log(-s) + \half +
  \order{\frac{1}{s}}  , 
\end{gather*}
we can write the \myrevision{intermediate-range} solution with terms
ordered by size (for this range of $s$) as
\begin{multline}
  v(s) = \left( \half-s \right) + 
  \frac{1}{\gamma^2}  \left(\half - s\right) \log(-s) +
  c \;  s \frac{\log\gamma}{\gamma^2}  + \frac{1}{\gamma^2} \left(d+\half\right)s\; +
  \\
  a \frac{\log\gamma}{\gamma^2} + \frac{b}{\gamma^2} +
  \frac{1}{\gamma^2} \half(1+\log 2) + \order{\frac{1}{\gamma^2
      s},\frac{1}{\gamma^4}}. 
  \label{eq:mid-ordered}
\end{multline}
Next, rewrite the inner solution \en{vinn} using $v = \gamma^2
\vinn(r)$ and 
\begin{gather*}
  r = 1+\frac{s}{\gamma^2}.
\end{gather*}
Then, for $|s| = \order{\gamma}$, we obtain the inner
expansion in terms of the intermediate range variable as
\begin{gather*}
  \gamma^2\vinn\left( 1+\frac{s}{\gamma^2} \right) = 
  \left(\half-s\right) - \frac{1}{\gamma^2} s\log(-s) +
  2s\frac{\log\gamma}{\gamma^2} + 
  \frac{17}{12 \gamma^2} + \frac{1}{\gamma^2} \half\log(-s) -
  \frac{\log\gamma}{\gamma^2} + \order{\frac{1}{\gamma^4}}.  
\end{gather*}
By comparing this expression with that for the intermediate range in
\en{mid-ordered}, we find that all terms match (both constant and
logarithmic) provided the constants satisfy the identities
\begin{gather*}
  a = -1, \qquad 
  b = \frac{11}{12} - \half \log 2, \qquad
  c = 2, \qquad
  d = -\half.
\end{gather*}
Therefore, the \myrevision{intermediate-range} solution for $s < 1/2$ is
given asymptotically by
\begin{gather}
  v(s) = \left(\half-s\right) + (2s-1)\frac{\log\gamma}{\gamma^2} +
  \frac{11}{12 \gamma^2} + \frac{1}{\gamma^2} \left(\half-s\right)
  \log\left(\half-s\right) + \order{\frac{1}{\gamma^4}},
  \label{eq:vmid-left}
\end{gather}
which we can write as
\begin{gather}
  v(s) = \left( \frac{1}{2} + \frac{11}{12 \gamma^2} - s \right)
  \left(1 + \frac{-2 \log(\gamma) + \log(1/2 - s)}{\gamma^2} \right) +
  \order{\frac{1}{\gamma^4}}. 
  \label{eq:allinner}
\end{gather}
Substituting the constants into \en{v0-mid}, we also find that to
leading order 
\begin{gather}
  v(0) = \frac{1}{2} + \frac{ b - \log(\gamma)}{\gamma^2} +
  \order{\frac{1}{\gamma^4}} \qquad \text{or} \qquad 
 u(1) = \frac{1}{2\gamma^2} +
  \frac{b-\log\gamma}{\gamma^4} + \order{\frac{1}{\gamma^6}}.
  \label{eq:vmid0-asy}
\end{gather}
Similarly, \en{vinf-mid} implies that 
\begin{gather*}
  \log\vinf = -\gamma^2 + 2 \log\gamma - \half +
  \order{\frac{1}{\gamma^2}}, 
\end{gather*}
so that
\begin{gather}
  \vinf = \gamma^2 e^{-\gamma^2-1/2 +\order{1/\gamma^2}}  
  \qquad \text{or} \qquad 
  \eps = \uinf = e^{-\gamma^2-1/2 +\order{1/\gamma^2}}.  
  \label{eq:vmidinf-asy}
\end{gather}

\subsubsection{A More Precise Estimate of Front Location and 
  Solution Curvature}
\label{sec:mid-mid}

The expression \en{allinner}, in which $v(s)$ would appear to vanish
when $s = 1/2 + 11/12\gamma^2 + \order{1/\gamma^4}$, is strongly
suggestive of a wetting front location given in terms of the inner
variable by
\begin{gather}
  \rstar  = 1 + \frac{1}{2 \gamma^2} + \frac{11}{12 \gamma^4} +
  \order{\frac{1}{\gamma^6}}. 
  \label{eq:wettingc}
\end{gather}
To further support this estimate, we return to the expression
\en{midmatch1} for the overall behaviour of the solution in the
intermediate layer, and continue the calculation given earlier using
the more refined values for $u(1)$ and $\uinf$ given by
\en{vmid0-asy} and \en{vmidinf-asy}.  Expanding the $\li(z)$ function
to include the $z/\log^2 z$ term, we obtain after some manipulation
\begin{align*}
  \gamma^2 \uinf \; \li \left( \frac{u(1)}{\uinf} \right) &= 
  \frac{\frac{1}{2} + \frac{1}{\gamma^2}\left( \frac{11}{12}  -
      \frac{\log(2)}{2} - \log(\gamma)\right)}{\log(1/2\gamma^2) +
    \gamma^2 + 1/2} 
  + \frac{1/2}{\gamma^4} + \order{\frac{1}{\gamma^6}},\\
  &=\frac{1}{2\gamma^2} + \frac{11}{12 \gamma^4} + \frac{1}{4\gamma^4} +
  \order{\frac{1}{\gamma^6}}.
\end{align*}
Thus, as
\begin{gather*}
  \rstar(r-1) = \gamma^2 \uinf \; \li \left(
    \frac{u(1)}{\uinf} \right) - \gamma^2 \uinf \; \li \left(
    \frac{u(r)}{\uinf} \right),
\end{gather*}
we have
\begin{align*}
  (r-1) &= \frac{\frac{1}{2\gamma^2} + \frac{11}{12 \gamma^4} +
    \frac{1}{4\gamma^4} + \order{\frac{1}{\gamma^6}}} {1 +
    \frac{1}{2\gamma^2} + \order{\frac{1}{\gamma^4}}}- \gamma^2
  \uinf \; \li \left( \frac{u(r)}{\uinf} \right)\\
  &= \frac{1}{2\gamma^2} + \frac{11}{12 \gamma^4} +
  \order{\frac{1}{\gamma^6}} - \gamma^2 \uinf \; \li \left(
    \frac{u(r)}{\uinf} \right).
\end{align*}
We deduce that if $\rstar$ is now as given in \en{wettingc}, then with
this refined value
\begin{gather}
  r-\rstar = - \gamma^2 \uinf \; \li \left( \frac{u(r)}{\uinf}
  \right). 
  \label{eq:inter12}
\end{gather}
Applying the previous reasoning, this leads to exactly the expression
\en{inter2}, with $\rstar$ now given by the refined approximation
\en{wettingc}. 

An alternative definition of the wetting front location is that it is
the point $\rstars$ where the curvature $\dderiv{u}{r}$ takes its maximum
value, so that $u^{\prime\prime\prime}(\rstars) = 0$. We now show that
this is equivalent to the value $\rstar$ given above, and estimate the
curvature at this point. Differentiating the underlying differential
equation yields
\begin{gather*}
  u\ddderiv{u}{r} + \deriv{u}{r} \dderiv{u}{r} =
  -r\dderiv{u}{r}/\gamma^2 - \deriv{u}{r}/\gamma^2,
\end{gather*}
so that imposing $u^{\prime\prime\prime}(\rstars)=0$ yields at $r =
\rstars$ 
\begin{gather*}
  \deriv{u}{r} = -
  \frac{r\dderiv{u}{r}}{1+ \gamma^2 \dderiv{u}{r}}. 
\end{gather*}
Assuming that $\rstars \approx 1$ and $\dderiv{u}{r}(\rstars)\gg 1$, we
have to leading order
\begin{gather*}
  \deriv{u}{r} (\rstars) = -\frac{1}{\gamma^2}.
\end{gather*}
However, based on the intermediate range equation derived earlier, we
also have 
\begin{gather*}
  \gamma^2 \deriv{u}{y} = -\log\left(\frac{u}{\uinf}\right),  
\end{gather*}
and so it follows immediately that the maximum value of $\dderiv{u}{r}$
arises when
\begin{gather}
  \log\left(\frac{u}{\uinf}\right) = 1 
  \qquad \text{or} \qquad
  u(\rstars) = e \uinf.
  \label{eq:a1}
\end{gather}
Substituting these various approximations into the underlying ODE we
find that to leading order the maximum value of $\dderiv{u}{r}$ at
$r=\rstars$ is given by
\begin{gather}
  \max \dderiv{u}{r} \approx \frac{1}{\gamma^4 e \uinf} =
  \frac{e^{\gamma^2-1/2}}{\gamma^2}.
  \label{eq:max-Theta-yy}
\end{gather}
As expected, this is very large indeed for even moderately large
$\gamma$. This incidentally leads to severe problems in any numerical
scheme that attempts to compute for large $\gamma$.

The wetting front location may now be estimated by substituting \en{a1}
into \en{inter12}  to give
\begin{gather*}
  \rstars - \rstar = - \gamma^2 \uinf \li(e). 
\end{gather*}
Because $\li(e)$ is of order one, we have that $\rstars$ and $\rstar$
are identical to all polynomial orders.

\subsection{Outer Solution, and Matching to the Intermediate Solution}
\label{sec:mid-outer}

To complete this calculation, we now match to the outer solution.  We
recall from \en{inter2} that the outer form of the intermediate range
solution is given by
\begin{gather}
  u = \uinf \exp \exp \left( -\frac{(r - \rstar)}{\gamma^2
      \uinf} - \eulergam \right) . 
  \label{eq:inter3}
\end{gather}
We observe from this expression that $u$ is very close to $\uinf$
if $r$ is only {\em slightly larger} than $\rstar$. In this range, the
term $\exp \left( -\gamma^{-2} e^{\gamma^2} \left( r - \rstar \right) -
  \eulergam \right)$ is very small, so that we may approximate
\en{inter3} by 
\begin{gather}
  u = \uinf \left[ 1 +  \exp \left( -
    \frac{(r-\rstar)}{\gamma^2 \uinf} - \eulergam \right) \right]. 
  \label{eq:inter4}
\end{gather}

Within the outer region, we have that $u \to \uinf$ as $r\to\infty$.
To match to the solution outer region, we take this to correspond to
those $r$ values for which $u$ is close to $\uinf$ and look for a
solution that is a small perturbation from a constant, namely
\begin{gather}
  u(r) = \uinf (1 + g(r)), \label{eq:thout1}
\end{gather}
where both $|g|$, $|\deriv{g}{r}| \ll 1$ for $y$ sufficiently large.
After substituting this expression into \en{usim-u}, we obtain
\begin{gather*}
  \gamma^{-2}  \deriv{g}{y} = -\frac{\uinf}{r}\, (1+ g)\, \dderiv{g}{r}, 
\end{gather*}
%
%
which can be approximated for small $|g|$ by 
\begin{gather*}
  \uinf \dderiv{g}{r} = - \gamma^{-2} r \deriv{g}{yr}.
\end{gather*}
This equation can be integrated once to obtain
\begin{gather}
  \deriv{g}{y} = -A \exp\left( -\frac{
      \left(r^2-\rstar^2\right)}{2 \uinf \gamma^{2}} \right),  
  \label{eq:wp-soln}
\end{gather}
where $A>0$ is a constant. Close to $\rstar$ we then have
\begin{gather}
  \deriv{g}{y} = -A \exp\left( -\frac{\rstar (r-\rstar)}{
      \uinf \gamma^{2} } \right).
  \label{eq:wp1-soln}
\end{gather}
Because $\rstar=1$ to leading order, this expression for $g'$ exactly
matches the derivative of the expression \en{inter4} as long as
\begin{gather}
  A = \frac{e^{-\eulergam}}{\uinf\gamma^{2}}.  
  \label{eq:A-constant}
\end{gather}
This completes the outer solution matching to leading order.

\subsection{Asymptotic Expansion for $\gamma$}
\label{sec:gamma}

We have so far expressed all asymptotic expansions in terms of the large
parameter $\gamma$.  However, $\gamma$ is not actually known
\emph{a~priori} for a physical problem and so instead it is preferable
to write $\gamma$ in terms of the known parameter $\bbar$.  Taking the
result from \en{vmidinf-asy}, we have that
\begin{gather}
  \bbar = \gamma^2 + \half + \frac{\alpha_3}{\gamma^2} +
  \order{\frac{1}{\gamma^4}},
  \label{eq:bbar-gamma}
\end{gather}
where we have introduced the coefficient $\alpha_3$ that is yet to be
determined.  Neglecting the $\order{\gamma^{-4}}$ terms yields a
quadratic equation for $\gamma^2$ whose solution can be expressed for
large $\bbar$ as
\begin{gather}
  \gamma = \bbar^{1/2} - \frac{1}{4}\bbar^{-1/2} -
  \left(\frac{\alpha_3}{2}+\frac{1}{32}\right) \bbar^{-3/2} +
  \order{\bbar^{-5/2}}.
  \label{eq:gamma-bbar}
\end{gather}
The high fidelity numerical calculations reported by Amodio
et~al.~\cite{amodio-etal-2014} provide convincing numerical evidence
that 
\begin{gather*}
  \alpha_3 = \frac{1}{12},
\end{gather*}
although we have no formal calculation as yet to confirm this result.


\subsection{Summary of Composite Asymptotic Solution}
\label{sec:matched-plot}

Here we present a concise summary of the leading order asymptotic
approximations for the inner, \myrevision{intermediate} and outer
solutions derived in Eqs.~\en{vinn}, \en{vmid-left}, \en{inter4}, 
\en{thout1}, \en{wp-soln} and \en{A-constant}, expressing each in terms of
the original similarity variables $y$ and $\Theta(y)$ as follows:
\begin{subequations}\label{eq:compos}
  \begin{alignat}{4}
    \begin{array}[t]{c}
      \text{Inner:}\\
      \text{($0\leqslant y<\gamma^{-1}$)} 
    \end{array}
    & \quad & \Theta(y) & = 
    1-\gamma y + \frac{1}{\gamma^2} \left[ \half - \half(1-\gamma y)^2
      + (1-\gamma y) \log(1-\gamma y) \right] + {\cal O}\left(\frac{1}{\gamma^4}\right),
    \label{eq:compos-in}\\
    \begin{array}[t]{c}
      \text{Left \myrevision{intermediate}:}\\
      \text{($\gamma^{-1} < y < \ystar$)}
    \end{array}
    & & \Theta(y) & =
    \frac{1}{\gamma} (\ystar-y) \, \Big[ \gamma^2 + \log \gamma +
        \log(\ystar-y) \Big] + \order{\frac{1}{\gamma^3}} 
    \label{eq:compos-midL},\\
    \begin{array}[t]{c}
      \text{Right \myrevision{intermediate}:}\\
      \text{($y > \ystar$)}
    \end{array}
    & & \Theta(y) & \approx 
    \eps + \eps \exp \left( -\eulergam - \frac{y-\ystar}{\gamma\eps}
      \right)
    \label{eq:compos-midR},\\
    \begin{array}[t]{c}
      \text{Outer:}\quad\mbox{}\\
      \text{($y \to \infty$)} \mbox{\quad}
    \end{array}
    & & \Theta'(y) & \approx
    \frac{1}{\gamma} \frac{e^{-\eulergam}}{\eps} \exp\left(
      - \frac{y^2-\ystar^2}{2\eps} \right)
    \label{eq:compos-out}.
  \end{alignat}
  In these above expressions, $\eps=e^{-\bbar}$ and $\gamma$ can be
  written in terms of $\bbar$ as 
  \begin{gather}
    \gamma \approx \bbar^{1/2} - \frac{1}{4}\,\bbar^{-1/2} 
    - \left(\frac{\alpha_3}{2}+\frac{1}{32}\right)\,\bbar^{-3/2} 
    . \label{eq:compos-zo-gamma}
  \end{gather}
  This last equation may also be used to replace $\gamma$ in the 
  expansion of the front location~\en{wettingc} to yield 
  \begin{gather}
    \ystar \approx \bbar^{-1/2} + \frac{3}{4} \bbar^{-3/2} +
    \left(\frac{\alpha_3}{2} + \frac{133}{96}\right) \bbar^{-5/2} 
    , \label{eq:compos-ystar-beta}
  \end{gather}
\end{subequations}
where $\alpha_3=\frac{1}{12}$ (numerically).

\section{Other Asymptotic Approximations}

In this section, we present two alternate asymptotic solutions to the
exponential diffusion problem that are derived using other methods.

\subsection{Babu's Solution}
\label{sec:babu}

Babu~\cite{babu-1976a} studied a modified version of the exponential
diffusion problem in which he assumed that beyond some distance
$y\geqslant\ystar$ the saturation is equal to the residual value.  In
other words, he solves a modified boundary value problem that replaces
\myrevision{our condition at infinity~\en{usim-bcL} with
\begin{gather}
  \Theta(\ystar) = \eps \qquad \text{and} \qquad
  \frac{d\Theta}{dy} \to 0 
  \quad \text{as $y\to\ystar$} 
  \tag{\ref{eq:usim-bcL}${}^{\prime\prime}$}
  \label{eq:ubabu-bcL}.
\end{gather}
Here, the front location $\ystar$ is an unknown constant that must be
determined as part of the solution, which explains the need for the
extra limiting condition in \en{ubabu-bcL} that corresponds physically
to requiring a zero water flux at $\ystar$.}  By expanding the solution
as a series in $\eta=(1-\eps)/\bbar$, Babu obtained the following
approximation for the wetting front location
\begin{subequations}\label{eq:babusolution}
  \begin{align}
    \ystar^B &= \eta^{1/2}
    \left[ 1 + \frac{1}{3}\eta
      + \left( \frac{17}{90} + \frac{\bbar}{8} \right) \eta^2
      + \cdots \right],
    \label{eq:ystar-babu}
  \end{align}
  where we use a superscript ``B'' to distinguish Babu's solution.
  Neglect exponentially small terms in $\eps=e^{-\bbar}$, this
  expression may be rewritten in terms of $\bbar$ as
  \begin{align}
    \ystar^B & \approx \bbar^{-1/2} + \frac{11}{24} \bbar^{-3/2}
    + \order{\bbar^{-5/2}} , 
    \label{eq:ystar-babu2}
  \end{align}
  which may then be compared directly to our asymptotic approximation in
  \en{compos-ystar-beta}. Although the leading order terms in
  $\order{\bbar^{-1/2}}$ are identical, the difference in coefficients
  at the next order generates a significant discrepancy in the front
  locations.  Babu also derived two further approximations that he
  refers to as second- and third-order expansions
  \begin{align}
    \Theta^{B,2}(y) &= 
    \displaystyle 1 + \bbar\eta^{1/2} \,
    \left( \frac{\eta}{6}\, y -y + \frac{1}{6}\,y^3
    \right),  
    \label{eq:theta-babu2}\\
    \Theta^{B,3}(y) &= 
    \Theta^{B,2}(y) + 
    \bbar\eta^{1/2} \, \left[  
      \left( \frac{7}{360}+\frac{\bbar}{24} \right) \eta^2 y
      - \frac{\eta}{36} y^3 + \frac{\bbar \eta^{1/2}}{12} y^4 
      - \frac{1}{40} y^5 
    \right],  
    \label{eq:theta-babu3}
  \end{align}
\end{subequations}
both of which are defined on the interval $0\leqslant y \leqslant
\ystar^B$. When $y>\ystar^B$, the saturation in both cases satisfies
$\Theta^{B,2} = \Theta^{B,3} = \eps$ (corresponding to $\sbar=0$).
Babu's solution is simpler than the one derived in this paper in that it
is a polynomial expansion in integer powers of $y$ and requires no
multi-layer matching.  However, this simplicity comes at the expense of
reduced accuracy as well as a lack of information about the detailed
structure of the wetting front. Another disadvantage of Babu's approach
is that none of his approximations for $\Theta$ is continuous at
$y=\ystar^B$.

\subsection{Parlange's Solution}
\label{sec:parlange}

\myrevision{Many approaches for solving the nonlinear diffusion equation
  analytically begin with the \emph{Bruce-Klute equation},
  $D(\theta) = - \frac{dy}{d\theta} \; \int_{\theta_o}^{\theta} y(\alpha)\,
  d\alpha,$
  which can be derived from~\en{usim} by treating $y$ as a function of
  saturation $\theta$~\cite{bruce-klute-1956}.  Although this equation
  was originally used to help interpret experimental data, it has
  subsequently been exploited by many authors to derive analytical
  solutions of the nonlinear diffusion equation, most notably by
  Philip~\cite{philip-1960, philip-1973}.  A related iterative solution
  approach was developed by Parlange~\cite{parlange-1971i,
    parlange-1973}}, who applied a number of simplifications to obtain
an approximate formula for $y$ in terms of integrals of the diffusivity.
For the specific case of an exponential diffusivity,
Parlange~\cite{parlange-1973} derived an implicit asymptotic
representation for the saturation that can be expressed in our
similarity variables as
\begin{subequations}\label{eq:parlange}
  \begin{gather}
    y = \left(\bbar^{-1/2}+\bbar^{-3/2}\right) \left( 1 - \Theta^P + \Theta^P
      \log \Theta^P / \bbar \right).
    \label{eq:parlange-theta}
  \end{gather}
  In a similar manner to Babu, Parlange also imposed the
  condition that $\Theta\equiv\eps=e^{-\bbar}$ for all
  $y\geqslant\ystar^P$.  Equating terms in Eq.~\en{parlange-theta} gives
  the result
  \begin{gather}
    \ystar^{P} =  \left(\bbar^{-1/2}+\bbar^{-3/2}\right) \left( 1 -
      2 e^{-\bbar}\right).
    \label{eq:parlange-ystar}
  \end{gather}
\end{subequations}
Intriguingly, this estimate involves an exponential term that would
correspond to a ``beyond all orders'' contribution to the results in
this paper. This estimate behaves asymptotically in the limit of large
$\bbar$ as
\begin{gather}
  \ystar^P \approx \bbar^{-1/2} +\bbar^{-3/2},
  \label{eq:cpar}
\end{gather}
from which it is clear that the front location matches only at leading
order with our estimate \en{compos-ystar-beta} and that of Babu from
\en{ystar-babu2}.

Another related solution was derived from the Bruce-Klute equation by
Parlange~\etal~\cite{parlange-etal-1992}, who extended an earlier method
of Heaslet and Alksne~\cite{heaslet-alksne-1961} for the power-law
diffusivity to the exponential case.  This solution also has
representation for $\Theta$ that involves the exponential integral
function, $\Ei(z)$.  Because evaluating $\Theta$ requires inverting
$\Ei$, their solution is much more complicated and expensive to compute
than ours and so we have not considered it here.

\section{Comparison of Asymptotic and Numerical Results}
\label{sec:numerical}

We now compare the various asymptotic solutions presented in the
preceding sections.  We also validate the asymptotic results using
numerical simulations of the initial value problem in
Eqs.~\enprime{usim} based on an algorithm that is described next.

\subsection{Solution Algorithm}
\label{sec:numerical-algo}
The ODE~\en{usim-u} for the similarity variable $\Theta(y)$ is solved
numerically over an interval \myrevision{$y\in[0,M]$} using the \Matlab\
initial value solver {\tt ode15s}.  Because both the problem and our
asymptotic results correspond \myrevision{to the semi-infinite interval
  $[0,\infty]$}, we must choose \myrevision{$M$} sufficiently large that
any error arising from truncating the right-hand boundary does not
pollute the solution in the interior; on the other hand, $\Theta(y)$
tends quite rapidly toward $\eps$ as $y$ increases beyond the front, and
so \myrevision{$M$} does not need not be taken much larger than
$\ystar$.  In practice, we have found that choosing
$\myrevision{M}=2\bbar^{-1/2}$ (twice the leading order estimate of the
wetting front location) provides a reasonable compromise between
efficiency and accuracy.

For an actual wetting scenario, we know the asymptotic saturation $\eps$
(from $\bbar$), but the value of $\gamma$ in the second initial
condition is not known \emph{a priori}.  We therefore build the initial
value solver into a shooting type algorithm, for which we guess the
value $\gamma=-\deriv{\Theta}{y}(0)$, integrate the saturation variable
to \myrevision{$\Theta(M)$}, and then compare to the target value
$\eps$.  The value of $\gamma$ is then modified using \myrevision{the
  bisection method} and this integration procedure is iterated until the
relative error in the right boundary condition satisfies
$|\myrevision{\Theta(M)}-\eps|/\eps<\tol$, where $\tol$ is a given
tolerance.  The leading term $\gamma\approx \bbar^{1/2}$ from the
asymptotic formula \en{compos-zo-gamma} provides a sufficiently accurate
initial guess to begin the iteration, and choosing ODE tolerances of
\abstol\ $=$\ \reltol\ $=10^{-10}$ and $\Theta$-tolerance of
$\tol=10^{-6}$ yields a computed solution that for all intents and
purposes can be treated as an ``exact solution'' of the original
problem.

A significant difficulty with this algorithmic approach is that for even
moderate values of $\gamma$ the solution curvature near the wetting
front is on the order of $\exp(\gamma^2)$, which can be extremely
large. This large curvature causes problems with the initial value
solver we are using, and effectively limits the maximum value of $\bbar$
to roughly 20, which restricts $\gamma \lessapprox 4.4$.  With this
restricted range of $\gamma$, we can still test the validity of the
asymptotic formulae derived here but we are not able to see their true
accuracy at larger $\gamma$, as the asymptotic series we have derived
only converge at a polynomial rate as $\gamma$ increases.  Consequently,
we will also report some results using a much more sophisticated
numerical approach that is the subject of a related
paper~\cite{amodio-etal-2014} and which employs a high-order Taylor
expansion based boundary value solver coupled with mesh adaptivity. This
method permits calculations up to $\gamma=18$, corresponding to
$\bbar\approx 325$ and $\theta_{\infty} \approx 1.18\times 10^{-141}$.


\subsection{Saturation Profiles}
\label{sec:numerical-sat}
Plots of the multi-layer asymptotic solution determined by
Eqs.~\en{compos} are depicted in Figure~\ref{fig:compare} for values of
$\bbar=4$, 8 and 16 alongside the corresponding plots of the computed
solution using the shooting method described above (which can
essentially be considered as an ``exact solution'').  The inner and
\myrevision{intermediate-range} solution are both displayed, and the
loss of asymptotic validity of the inner solution due to the logarithmic
term is evident as $y\to \ystar$.
We have included in all plots the corresponding asymptotic solutions of
Babu~\cite{babu-1976a} and Parlange~\cite{parlange-1973}, and in each
case a second plot of all curves in terms of the rescaled saturation
$\sbar$ is given (lower plots), which accentuates differences in the
wetting front location.  Our asymptotic solution is clearly an
improvement over Babu's solution for all values of $\bbar$, which we
attribute in large part to the fact that Babu's approximation truncates
the saturation at some approximate wetting front location and ignores
the details of the solution structure within and to the right of the
front.
\newcommand{\solidline}{\hdashrule[0.5ex]{1.2cm}{0.04cm}{}}
\newcommand{\dashedline}{\hdashrule[0.5ex]{1.3cm}{0.04cm}{0.2cm}}
\newcommand{\dashdotline}{\hdashrule[0.5ex]{1.2cm}{0.04cm}{0.15cm 0.05cm 0.05cm 0.05cm}}
\begin{figure}[tbhp]
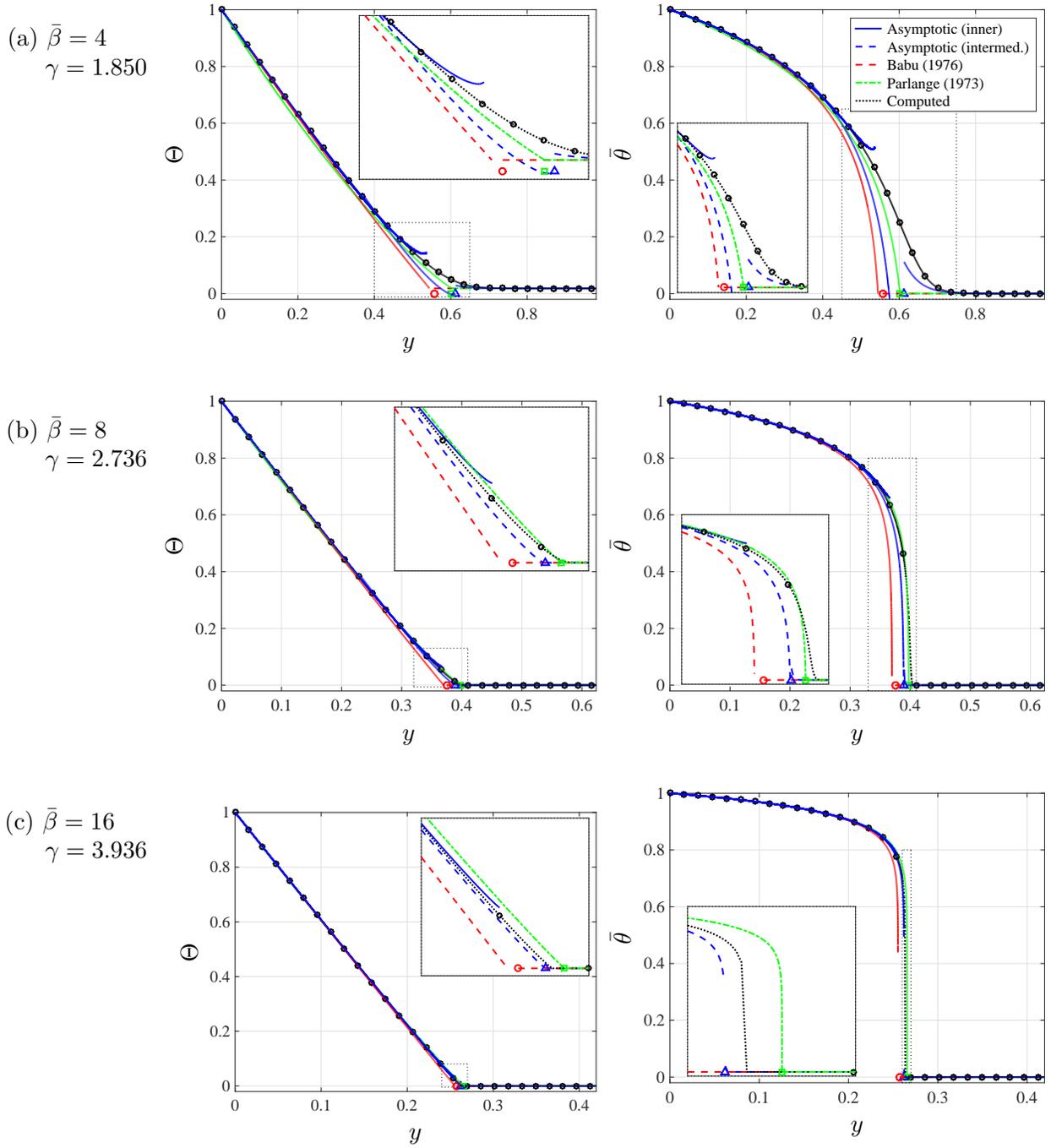

  \centering
  \psfrag{y}[cc][bc]{$y$}
  \psfrag{thetabar}[bc][cc]{$\sbar$}
  \psfrag{Theta}[bc][cc]{$\Theta$}
  \begin{tabular}{>{\flushright}p{0.14\textwidth}cc}
    (a) 
    $\bbar=4$\ \ \quad\break
    $\gamma=1.850$ 
    & 
    \raisebox{-\height}{\includegraphics[width=0.40\textwidth,clip]{compare4zoom.eps}} & 
    \raisebox{-\height}{\includegraphics[width=0.40\textwidth,clip]{compare4barzoom_MOD.eps}}
    \\
    & & \\
    (b)
    $\bbar=8$\ \ \quad\break
    $\gamma=2.736$ 
    &
    \raisebox{-\height}{\includegraphics[width=0.40\textwidth,clip]{compare8zoom.eps}} &
    \raisebox{-\height}{\includegraphics[width=0.40\textwidth,clip]{compare8barzoom.eps}}
    \\
    & & \\
    (c)
    $\bbar=16$\ \quad\break 
    $\gamma=3.936$  
    & 
    \raisebox{-\height}{\includegraphics[width=0.40\textwidth,clip]{compare16zoom.eps}} & 
    \raisebox{-\height}{\includegraphics[width=0.40\textwidth,clip]{compare16barzoom.eps}}
  \end{tabular}
  \caption{Comparison of our asymptotic solution
    (\textcolor{blue}{\solidline}) to those of Babu
    (\textcolor{red}{\dashedline}) and Parlange
    (\textcolor{green}{\dashdotline}) for $\bbar=4$, 8 and 16, displayed
    in terms of the similarity variable (left) and rescaled saturation
    (right).  The corresponding approximations of the wetting front
    location are denoted by points lying on the $y$-axis for our
    asymptotics (\textcolor{blue}{\pmb{\small$\triangle$}}), Babu
    (\textcolor{red}{\pmb{\LARGE$\circ$}}) and Parlange
    (\textcolor{green}{\pmb{\small$\square$}}).
  }
  \label{fig:compare}
\end{figure}

As described above, the relatively slow convergence of our asymptotic
solution, which is polynomial in $\gamma$, and hence also in $\bbar$,
means that our asymptotic approximation is not particularly accurate for
smaller values of $\bbar$ such as $\bbar=4$. This is reflected in the
mismatch between the left and right \myrevision{intermediate-range}
solutions at the wetting front location (represented by the blue
triangular point).  However, our asymptotic approximation improves
significantly in accuracy as $\bbar$ increases, and when we take
$\bbar=16$ the comparison between the computed and asymptotic results
are good.  The differences at even larger $\bbar$ are difficult to
visualize using saturation plots, and so we compare the solutions
further in the next sections in terms of the estimates for wetting front
location.

\subsection{Wetting Front Location, $\ystar$}
\label{sec:numerical-ystar}
\begin{figure}[tbhp]
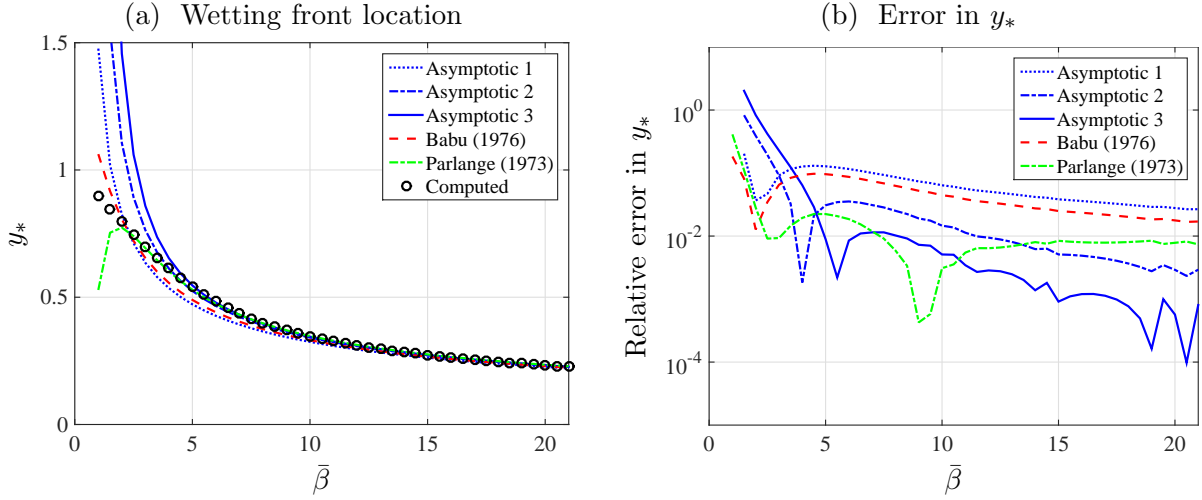

  \centering
  \psfrag{bbar}[cc][bc]{$\bbar$}
  \psfrag{ystar}[bc][cc]{$\ystar$}
  \psfrag{Rel. error in ystar}[bc][cc]{Relative error in $\ystar$}
  \begin{tabular}{ccc}
    (a)\ \ Wetting front location
    & & (b)\ \ Error in $\ystar$ \\
    \includegraphics[width=0.45\textwidth,clip]{yscompare3}
    & \quad & 
    \includegraphics[width=0.45\textwidth,clip]{yserror3}
  \end{tabular}
  \caption{Comparison of the various asymptotic estimates for wetting 
    front location $\ystar$ to the computed (``exact'') solution.}
  \label{fig:ystarcompare}
\end{figure}

We next focus on calculations of the wetting front location $\ystar$,
which we have defined in our derivation to the be point where
$\dderiv{\Theta}{y}$ is a maximum.
A visual comparison is provided in Figure~\ref{fig:ystarcompare} in
terms of the computed front location as a function of $\bbar$, showing
our asymptotic estimate $\ystar^{\text{\emph{asy}}}$ from
\en{compos-ystar-beta}, Babu's estimate $\ystar^{B}$ given in
\en{theta-babu3}, and Parlange's estimate $\ystar^P$ given in \en{cpar}.
We also included numerical estimates of front location using the Matlab
solver {\tt ode15s}, for which we were able to compute up to a maximum
of $\bbar\approx 20$ before the ODE solver failed.  The plot of error in
front location (relative to the computed solution) clearly shows that
although Parlange's estimate is better than our asymptotic solution for
some values of $\bbar \gtrapprox 10$, our result is consistently
superior for larger $\bbar$; indeed, even our two-term asymptotic
solution surpasses Parlange's result when $\bbar\gtrapprox 13$.


%
%
Although our simple numerical approach fails for $\bbar > 20$, a more
sophisticated numerical algorithm has been implemented by Amodio
\etal~\cite{amodio-etal-2014} that yields accurate solutions for values
of $\bbar$ much larger.  Their results for a much wider range of $\bbar$
are summarized in Table~\ref{tab:amodio}, along with the corresponding
asymptotic results for $\ystar^{\text{\emph{asy}}}$, $\ystar^B$ and
$\ystar^P$.  The superior accuracy of our asymptotic solution for large
$\bbar$ when compared with Babu's and Parlange's approximations is
evident upon comparing the computed value of $\ystar$ to the three
asymptotic approximations.
\begin{table}[tbhp]
  \centering
  \caption{Comparison of our asymptotic results to computations using
    the method of Amodio et al.~\cite{amodio-etal-2014}.}
  \label{tab:amodio}
  \renewcommand{\extrarowheight}{0mm}
  \begin{tabular}{c|ccc|ccc}\hline
    & \multicolumn{3}{c|}{Computed results} & 
    \multicolumn{3}{c}{Asymptotic results}\\
    $\gamma$        & $\bbar$    & $\theta_{\infty}$ & $\ystar$ &
    $\ystar^{\text{\emph{asy}}}$ & $\ystar^B$ & $\ystar^P$  \\ \hline 
    2  & 4.559435  & 1.046797e--2  & 0.571747  & 0.54536   & 0.51539   & 0.57103 \\
    6  & 36.50238  & 1.403505e--16 & 0.16911   & 0.1689165 & 0.16759   & 0.170050\\
    10 & 100.5008  & 2.254440e--44 & 0.1005094 & 0.1004949 & 0.100205  & 0.100743\\
    18 & 324.5002  & 1.178490e--141& 0.0556417 & 0.0556410 & 0.0055591 & 0.055683\\\hline
  \end{tabular}  
\end{table}

\subsection{Initial Slope, \myrevision{$-\gamma$}}
\label{sec:numerical-gamma}
\begin{figure}[tbhp]
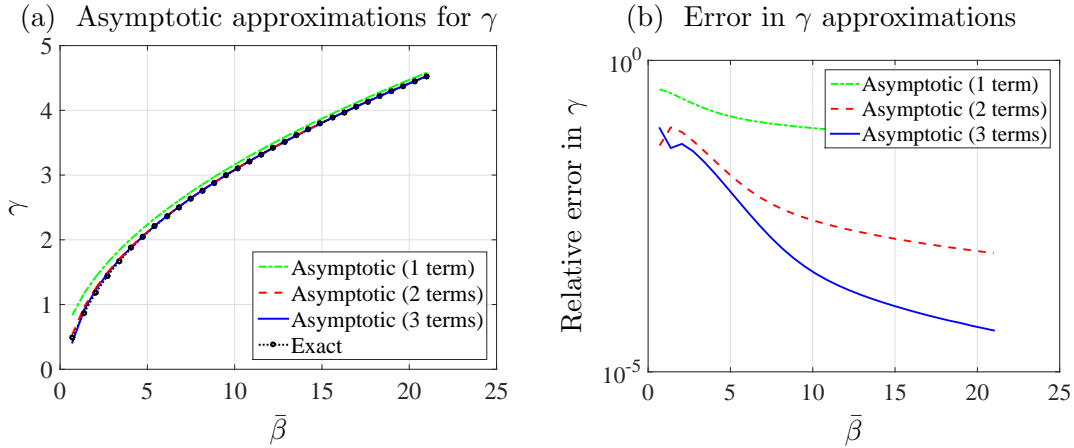

  \centering
  \psfrag{bbar}[cc][bc]{$\bbar$}
  \psfrag{gamma}[bc][cc]{$\gamma$}
  \psfrag{error}[bc][cc]{Relative error in $\gamma$}
  \psfrag{Error in gamma approximations}[cl][bl]{}
  \begin{tabular}{ccc}
    (a)\ \ Asymptotic approximations for $\gamma$ & & 
    (b)\ \ Error in $\gamma$ approximations\\
    \includegraphics[width=0.40\textwidth,clip]{shootgam.eps} 
    & \quad & 
    \includegraphics[width=0.40\textwidth,clip]{shootgamerr.eps}
  \end{tabular}
  \caption{Comparison of the various asymptotic approximations of
    $\gamma$ based on retaining one, two or three terms from the series
    in {\protect Eq.~\en{compos-zo-gamma}}.  The left hand plot compares
    the three asymptotic results to the value of $\gamma$ obtained using
    the shooting algorithm.  The right hand plot depicts the
    corresponding absolute error in the asymptotic approximations,
    treating the shooting results as the ``exact'' solution.}
  \label{fig:gamerr}
\end{figure}

Finally, we investigate the accuracy of our asymptotic approximation for
$\gamma$ in Eq.~\en{compos-zo-gamma}, \myrevision{where $-\gamma$ is the
  initial slope in the similarity variable $\Theta$.  No corresponding
  estimate is available from either Babu's or Parlange's solutions.}
For a range of $\bbar$, Figure~\ref{fig:gamerr} compares the value of
$\gamma$ determined from shooting simulations with that obtained from
the 1-, 2- and 3-term asymptotic estimates.  The accuracy of our series
approximation is evident for increasing values of $\bbar$, which is
essential because $\gamma$ is a required input for our numerical
simulations and yet it is the parameter $\bbar$ (not $\gamma$) that is
known \emph{a priori} for any given physical wetting scenario.

We finish by restating some of the numerical results obtained in
\cite{amodio-etal-2014} for values of $\gamma$ as high as 18
(corresponding to $\beta$ up to 324). Beyond this range, even this more
sensitive algorithm was unable to converge owing to the size of
$\theta''$ at the wetting front. The primary calculations of $\ystar$
and $\theta_{\infty}$ over a range of values of $\gamma$ were then
enhanced by using Richardson extrapolation. The primary contribution of
\cite{amodio-etal-2014} was to provide clear numerical evidence to
support the validity of the following expressions
\begin{gather*}
  \ystar \approx \frac{1}{\gamma} + \frac{1}{2 \gamma^3} + \frac{11}{12
    \gamma^5} + \frac{2.96}{\gamma^7}, 
\end{gather*}
and
\begin{gather*}
  \bbar = -\log (\theta_{\infty} ) \approx \gamma^2 + \frac{1}{2} +
  \frac{1}{12 \gamma^2} + \frac{0.089}{\gamma^4}, 
\end{gather*}
which extends the asymptotics derived in this paper by one additional
term.  These results are in full agreement with our asymptotics and also
hint at a more refined asymptotic calculation. 

\section{Conclusions}
\label{sec:conclude}

In this paper we have performed a multi-layer asymptotic analysis of a
nonlinear diffusion equation where the diffusion coefficient is an
exponential function, $D(\theta)=\Dzero e^{\beta\theta}$.
\myrevision{We focus on an application to the study of wetting front
  formation in unsaturated porous media flow, although exponential
  diffusion also arises in a number of other problems in heat transport,
  optical lithography, polymer diffusion, etc.}
\myrevision{The original boundary value problem for liquid saturation is
  reformulated as an initial value problem which, although not strictly
  necessary (c.f. Babu's solution~\cite{babu-1976a} in
  Section~\ref{sec:babu}), nonetheless permits an easy comparison to
  numerics using an initial value solver.}  Motivated by the structure
of the wetting front for large values of $\beta$, we used the method of
matched asymptotics to derive a four-layer solution consisting of a
different asymptotic series for each of four regions corresponding to
the wetting, and the leading and trailing edges.  Other previous
approaches to deriving (approximate) analytical solutions have assumed
that the reduced saturation is identically equal to zero at a finite
wetting front location, which essentially ignores the structure of the
{\itshape sharp corner} that appears at the front for large $\beta$.  In
contrast, our asymptotic solution uncovers the detailed structure of the
transition within the wetting front, where the solution has very high
curvature, and is considerably more accurate that other previous
approaches for large $\beta$.  The asymptotic solution converges
polynomially in $\beta$ and maintains a high degree of accuracy over a
wide range of $\beta$ values corresponding to physical porous media such
as soils and rock, unlike many other approaches which are either more
limited in their applicability or exhibit a logarithmic singularity near
the wetting front as $\beta$ becomes large.  Our approach has the
additional advantage that it expresses the solution in terms of explicit
formulas instead of requiring numerical approximation of a differential
equation (which is required in the approaches of
Wagner~\cite{wagner-1950} or Shampine~\cite{shampine-1973c}) or
inverting the exponential integral function (as in Parlange \etal's
asymptotic solution~\cite{parlange-etal-1992}).  \myrevision{In
  particular, our estimate for $\ystar$ may be of value to groundwater
  researchers and experimentalists since it provides an explicit formula
  for a measurable quantity (wetting front location) in terms of
  physical parameter values.}

\section*{Acknowledgments}
\label{sec:ack}

This work was funded by grants from the Natural Sciences and Engineering
Research Council of Canada, Mitacs Network of Centres of Excellence,
Pacific Institute for the Mathematical Sciences, and Engineering and
Physical Sciences Research Council.



\begin{thebibliography}{10}
\providecommand{\url}[1]{\texttt{#1}}
\providecommand{\urlprefix}{URL }
\expandafter\ifx\csname urlstyle\endcsname\relax
  \providecommand{\doi}[1]{doi:\discretionary{}{}{}#1}\else
  \providecommand{\doi}{doi:\discretionary{}{}{}\begingroup
  \urlstyle{rm}\Url}\fi
\providecommand{\eprint}[2][]{\url{#2}}

\bibitem{brutsaert-1979}
\textsc{W.~Brutsaert}, Universal constants for scaling the exponential soil
  water diffusivity?, \emph{Water Resour. Res.} 15:481--483 (1979).

\bibitem{clothier-white-1981}
\textsc{B.~E. Clothier} and \textsc{I.~White}, Measurement of sorptivity and
  soil water diffusivity in the field, \emph{Soil Sci. Soc. Amer. J.}
  45:241--245 (1981).

\bibitem{miller-bresler-1977}
\textsc{R.~D. Miller} and \textsc{E.~Bresler}, A quick method for estimating
  soil water diffusivity functions, \emph{Soil Sci. Soc. Amer. J.}
  41:1020--1022 (1977).

\bibitem{reichardt-nielsen-biggar-1972}
\textsc{K.~D. Reichardt}, \textsc{D.~R. Nielsen}, and \textsc{J.~W. Biggar},
  Scaling of horizontal infiltration into homogeneous soils, \emph{Soil Sci.
  Soc. Amer. Proc.} 36:241--245 (1972).

\bibitem{simpson-1993}
\textsc{W.~T. Simpson}, Determination and use of moisture diffusion coefficient
  to characterize drying of northern red oak (\emph{{Q}uercus rubra}),
  \emph{Wood Sci. Tech.} 27:409--420 (1993).

\bibitem{pel-1995}
\textsc{L.~Pel}, \emph{Moisture transport in porous building materials}, Ph.D.
  thesis Technische Universiteit Eindhoven (1995).

\bibitem{leech-lockington-dux-2003}
\textsc{C.~Leech}, \textsc{D.~Lockington}, and \textsc{P.~Dux}, Unsaturated
  diffusivity functions for concrete derived from {NMR} images,
  \emph{Mat\'{e}r. Constr.} 36:413--418 (2003).

\bibitem{cooper-1971}
\textsc{L.~Y. Cooper}, Constant temperature at the surface of an initially
  uniform temperature, variable conductivity half space, \emph{J. Heat
  Transfer} 93:55--60 (1971).

\bibitem{wagner-1950}
\textsc{C.~Wagner}, Diffusion of lead chloride in solid silver chloride,
  \emph{J. Chem. Phys.} 18:1227--1230 (1950).

\bibitem{hansen-1980}
\textsc{C.~M. Hansen}, Diffusion in polymers, \emph{Polymer Eng. Sci.}
  20:252--258 (1980).

\bibitem{wei-wuensch-1976}
\textsc{G.~C.~T. Wei} and \textsc{B.~J. Wuensch}, Tracer concentration
  gradients for diffusion coefficients exponentially dependent on
  concentration, \emph{J. Amer. Ceramic Soc.} 59:295--299 (1976).

\bibitem{vazquez-2007}
\textsc{J.~L. V{\'{a}}zquez}, \emph{The Porous Medium Equation: Mathematical
  Theory}, \emph{Oxford Mathematical Monographs} Clarendon Press, Oxford, 2007.

\bibitem{shampine-1973b}
\textsc{L.~F. Shampine}, Concentration-dependent diffusion. {II}. {S}ingular
  problems, \emph{Quart. Appl. Math.} 31:287--293 (1973).

\bibitem{crank-1975}
\textsc{J.~Crank}, \emph{The Mathematics of Diffusion}, second ed., Oxford
  University Press, 1975.

\bibitem{parslow-lockington-parlange-1988}
\textsc{J.~Parslow}, \textsc{D.~Lockington}, and \textsc{J.-Y. Parlange}, A new
  perturbation expansion for horizontal infiltration and sorptivity estimates,
  \emph{Transp. Porous Media} 3:133--144 (1988).

\bibitem{babu-1976a}
\textsc{D.~K. Babu}, Infiltration analysis and perturbation methods. 1.
  {A}bsorption with exponential diffusivity, \emph{Water Resour. Res.}
  12:89--93 (1976).

\bibitem{parlange-1973}
\textsc{J.-Y. Parlange}, A note on a three-parameter soil-water diffusivity
  function -- {A}pplication to the horizontal infiltration of water, \emph{Soil
  Sci. Soc. Amer. J.} 37:318--319 (1973).

\bibitem{parlange-etal-1992}
\textsc{M.~B. Parlange}, \textsc{S.~N. Prasad}, \textsc{J.-Y. Parlange}, and
  \textsc{M.~J.~M. R\"omkens}, Extension of the {H}easlet-{A}lksne technique to
  arbitrary soil water diffusivities, \emph{Water Resour. Res.} 28:2793--2797
  (1992).

\bibitem{elliott-etal-1986}
\textsc{C.~M. Elliott}, \textsc{M.~A. Herrero}, \textsc{J.~R. King}, and
  \textsc{J.~R. Ockendon}, The mesa problem: {D}iffusion patterns for $u_t =
  \nabla \cdot (u^m \nabla u)$ as $m\rightarrow +\infty$, \emph{IMA J. Appl.
  Math.} 37:147--154 (1986).

\bibitem{king-1988}
\textsc{J.~R. King}, Approximate solutions to a nonlinear diffusion equation,
  \emph{J. Eng. Math.} 22:53--72 (1988).

\bibitem{king-please-1986}
\textsc{J.~R. King} and \textsc{C.~P. Please}, Diffusion of dopant in
  crystalline silicon: {A}n asymptotic analysis, \emph{IMA J. Appl. Math.}
  37:185--197 (1986).

\bibitem{bear-1988}
\textsc{J.~Bear}, \emph{Dynamics of Fluids in Porous Media}, Dover, New York,
  1988.

\bibitem{yeh-franzini-1968}
\textsc{W.~W.-G. Yeh} and \textsc{J.~B. Franzini}, Moisture movement in a
  horizontal soil column under the influence of an applied pressure, \emph{J.
  Geophys. Res.} 71:5151--5157 (1968).

\bibitem{amodio-etal-2014}
\textsc{P.~Amodio}, \textsc{C.~J. Budd}, \textsc{O.~Koch},
  \textsc{G.~Settanni}, and \textsc{E.~B. Weinm\"uller}, Asymptotical
  computations for a model of flow in saturated porous media, \emph{Appl. Math.
  Comput.} 237:155--167 (2014).

\bibitem{bruce-klute-1956}
\textsc{R.~R. Bruce} and \textsc{A.~Klute}, The measurement of soil moisture
  diffusivity, \emph{Soil Sci. Soc. Amer. J.} 20:458--462 (1956).

\bibitem{philip-1960}
\textsc{J.~R. Philip}, General method of exact solution of the
  concentration-dependent diffusion equation, \emph{Aust. J. Phys.} 13:1--12
  (1960).

\bibitem{philip-1973}
\textsc{J.~R. Philip}, On solving the unsaturated flow equation: 1.~{T}he
  flux-concentration relation, \emph{Soil Sci.} 116:328--335 (1973).

\bibitem{parlange-1971i}
\textsc{J.-Y. Parlange}, Theory of water-movement in soils: {1}.
  {O}ne-dimensional absorption, \emph{Soil Sci.} 111:134--137 (1971).

\bibitem{heaslet-alksne-1961}
\textsc{M.~A. Heaslet} and \textsc{A.~Alksne}, Diffusion from a fixed surface
  with a concentration-dependent coefficient, \emph{J. Soc. Ind. Appl. Math.}
  9:584--596 (1961).

\bibitem{shampine-1973c}
\textsc{L.~F. Shampine}, Some singular concentration dependent diffusion
  problems, \emph{Z. angew. Math. Mech.} 53:421--422 (1973).

\end{thebibliography}

\providecommand{\noopsort}[1]{}

\end{document}